\documentclass[11pt]{article}
\usepackage[margin=1in]{geometry}
\usepackage{amsthm,amsmath,amsfonts,amssymb}
\usepackage[numbers,square]{natbib}
\usepackage{mathtools,booktabs}
\usepackage[colorlinks,citecolor=blue,urlcolor=blue]{hyperref}
\theoremstyle{plain}
\newtheorem{theorem}{Theorem}
\newtheorem{lemma}[theorem]{Lemma}
\newtheorem{corollary}[theorem]{Corollary}
\newtheorem{proposition}[theorem]{Proposition}
\newcommand{\E}{\mathbb E}
\newcommand{\Prob}{\mathbb P}
\newcommand{\CN}{\mathcal{CN}}
\newcommand{\diag}{\operatorname{diag}}
\newcommand{\tr}{\operatorname{tr}}
\newcommand{\op}{\mathrm{op}}
\newcommand{\ind}{\mathbf 1}
\hypersetup{pdftitle={Frieze's zeta(3) and the Avram--Bertsimas constant under Haar conjugation},pdfauthor={Jan Vecer}}
\title{Frieze's $\zeta(3)$ and the Avram--Bertsimas constant under Haar conjugation}
\author{Jan Vecer\\[0.25em]
  {\small Department of Probability and Mathematical Statistics, Charles University}\\
  {\small \href{mailto:Jan.Vecer@matfyz.cuni.cz}{Jan.Vecer@matfyz.cuni.cz}}}
\date{}

\begin{document}
\maketitle
\pagestyle{myheadings}
\markright{Minimum spanning trees under Haar conjugation}
\begin{abstract}
Let $H_d=U_d^*\Lambda_dU_d$, where $U_d$ is Haar unitary and
$\Lambda_d$ is a deterministic real diagonal matrix, and put
$\sigma_d^2=\mathbb E|H_{12}|^2$. For minimum spanning trees with
edge weights $\exp(-s|H_{ij}|^2/\sigma_d^2)$, we prove that the
normalized expected weights differ from those of the independent
Gaussian model by a polynomially small error, uniformly for
$0<s\le S$ with fixed $S>0$.
The limit is the classical Avram--Bertsimas constant and includes
Frieze's $\zeta(3)$ at $s=1$. For edge weights $|H_{ij}|^{-1}$,
we approximate the normalized expected weight by an explicit
one-dimensional integral with polynomial error, and obtain a divergent
expansion in inverse powers of $\log d$.
The assumptions require the spectrum, after centering and division
by its standard deviation, to remain uniformly bounded.
From Penrose's high-dimensional edge-length limit we derive
$2-\gamma-\zeta(2)$ as the first correction to the
high-dimensional Euclidean minimum spanning tree constant after
normalization by the unit-ball volume. This correction is a logarithmic
moment of the limiting distribution of MST edge weights; half of it is the first logarithmic correction
for the reciprocal Haar weights. Here $\gamma$ is Euler's constant
and $\zeta$ is the Riemann zeta function.
\end{abstract}

\par\smallskip
\noindent\textbf{MSC2020 subject classifications.} Primary 60C05; secondary 05C80, 60B20, 05C05.

\section{Introduction and main results}
We prove that classical minimum spanning tree constants for independent
edges remain valid under the dependence generated by Haar conjugation. Our main theorem includes Frieze's $\zeta(3)$
law and the entire family of
Avram--Bertsimas constants, with a
quantitative error uniform as the exponent tends to zero. This
uniformity gives asymptotics for reciprocal edge weights and connects
their first correction to high-dimensional Euclidean minimum spanning trees.

\citet{Frieze} proved that the expected minimum spanning tree (MST)
weight on the complete graph with independent uniform edge lengths
converges to $\zeta(3)$. \citet{Steele} extended the corresponding
convergence in probability beyond the original moment assumptions.
For edge distributions that behave like a power near zero,
\citet[Theorem~8]{AvramBertsimas} express the leading expected weight
through a family of constants. The independent-model series was also
obtained by \citet{Timofeev}, as they note explicitly. We prove that
the same constants describe expected MST weights under Haar conjugation,
with a quantitative error bound.

Our contributions are as follows.
\begin{enumerate}
\renewcommand{\labelenumi}{(\roman{enumi})}
\item \textbf{Expected MST weights under Haar dependence.}
We prove universality for weights $Q_e^s=\exp(-s|H_e|^2/\sigma_d^2)$
for every fixed $s>0$, including
Frieze's $\zeta(3)$ limit, with an explicit polynomial rate
(Theorem~\ref{thm:power}). Relative estimates for small component
events control component counts at moderate thresholds. A concentration argument over all
cuts supplies a global connectivity bound and makes those
estimates integrable over all edge thresholds.
\item \textbf{The reciprocal weights.}
We give an exact Laplace representation of the mean for reciprocal edge weights and derive an explicit integral profile with polynomial
error (Theorem~\ref{thm:universality}). We determine its asymptotic
coefficients and prove that the inverse-logarithmic series
diverges (Corollary~\ref{thm:expansion}). We also characterize the
spectral growth for which the MST weight vanishes. For the spectrum
$\lambda_k=k^{3/2}$, its expected weight multiplied by $\sqrt{\log d}$
tends to $10/3$.
\item \textbf{An explicit consequence of Penrose's limit.}
In spatial dimension $D$, we evaluate $2-\gamma-\zeta(2)$ as the
first $1/D$ coefficient after multiplying the Euclidean MST constant
by the $D$th root of the unit-ball volume
(Proposition~\ref{cor:geometric-endpoint}), using the high-dimensional
limit and tail estimates of \citet{PenroseHighDim}. Half this same
coefficient gives the first logarithmic correction in the Haar
model with reciprocal edge weights.
\end{enumerate}

The geometric result is a logarithmic-moment consequence of
\citet[Proposition~3 and Theorem~2]{PenroseHighDim}, who identifies
the limiting empirical distribution of Euclidean MST edge lengths
in high dimensions. His Lemmas~3--4 also provide uniform tail
estimates sufficient for the moment convergence. We evaluate that
moment in closed form and identify its relation to the Haar correction.
Our short component-count proof uses the earlier continuum-percolation
limit of \citet[Proposition~1 and Theorem~1]{Penrose} and an elementary
degree bound. The leading Euclidean equivalent was already established
by \citet{BertsimasVanRyzin}.

The proofs use the component-count integral of \citet{AvramBertsimas}.
\citet{Janson} developed a functional limit theorem for tree-component
counts to establish Gaussian fluctuations of the independent
uniform-weight MST; exact expectation formulas for general independent
edge laws are available from \citet{LiZhang}. Related expectation
identities and extensions are surveyed by \citet{SteeleMST} and
developed by \citet{FillSteele}. Our estimates address rare component events under
Haar dependence and the global tail needed to integrate them.

Recent work of \citet{CMcontrol,CM} provides a motivating special case:
the MST with reciprocal weights $|H_{ij}|^{-1}$ for the spectrum
$\lambda_k=k^{3/2}$. These weights arise in quantum control.
Under the structural hypotheses of \citet[Proposition~5]{ChambrionFlux},
when the consecutive eigenvalue gaps of the control operator are
rationally independent, the worst-case state-transfer
time is at most $\pi M_d/2$, where $M_d$ is the MST weight for these
reciprocal weights. \citet[Section~1.3]{CM} conjecture that
$M_d\to0$ in probability and that
$\E M_d\sim\alpha(\log d)^{-1/2}$ for some $\alpha>0$;
the surrogate calculation of \citet[Section~2.2]{CM}
yields $\alpha=10/3$.
Theorem~\ref{thm:main} and \eqref{eq:power-leading}--\eqref{eq:critical-second}
prove both conjectures and identify this constant. Their reciprocal-weight
problem is one application of the present comparison, which also covers
the family $Q_e^s$, including Frieze's constant, for general spectra
satisfying \eqref{eq:bounded}.
Equation~\eqref{eq:universality} gives a quantitative asymptotic
counterpart of their observed agreement between Haar and independent
means. \citet[Section~2.2, equation~(6)]{CM} also predict heuristically
a first correction near $-0.11$ for the independent model; we evaluate
it in closed form and prove it for both models.

Section~\ref{sec:spectral-growth} applies the asymptotics for reciprocal edge weights to regularly varying spectra.

\subsection{The models and edge weights}
For each integer $d\ge2$, let
$\Lambda_d=\diag(\lambda_{1,d},\ldots,\lambda_{d,d})$
be real diagonal and let $U_d$ be Haar distributed on the complex
unitary group $U(d)$. Set
\[
 H_d=U_d^*\Lambda_dU_d,\qquad
 \bar\lambda_d=\frac1d\sum_k\lambda_{k,d},\qquad
 v_d=\frac1d\sum_k(\lambda_{k,d}-\bar\lambda_d)^2,
 \qquad \sigma_d^2=\frac{dv_d}{d^2-1}.
\]
The last quantity is the exact off-diagonal second moment. Write
$H_{ij}=(H_d)_{ij}$. On the complete graph with vertex set
$\{1,\ldots,d\}$, give the edge $e=\{i,j\}$, $i<j$, weight
$W_e=W_{ij}=|H_{ij}|^{-1}$ and denote the minimum spanning tree
weight by $M_d$. We write $H_e=H_{ij}$ and use the same edge-index
convention for other arrays. All minima over $T$ below range over
spanning trees of this complete graph.

In the independent and identically distributed (iid) model, write
\[
 H_{ij}^{\mathrm{iid}}
 =\frac{\sigma_d}{\sqrt2}\bigl(N_{ij}^{(1)}+iN_{ij}^{(2)}\bigr),
 \qquad i<j,
\]
where all $N_{ij}^{(1)},N_{ij}^{(2)}$ are independent real $N(0,1)$
variables. This is the circular complex Gaussian law, denoted by
$\CN(0,\sigma_d^2)$; thus
$\E|H_{ij}^{\mathrm{iid}}|^2=\sigma_d^2$, and
$|H_{ij}^{\mathrm{iid}}|^2/\sigma_d^2$ is a mean-one exponential.
Set $H_{ji}^{\mathrm{iid}}=\overline{H_{ij}^{\mathrm{iid}}}$;
the diagonal is irrelevant. Denote its MST weight by
$M_d^{\mathrm{iid}}$. More generally, $\CN(0,\Sigma)$ denotes
the circular complex Gaussian vector law with covariance
$\E ZZ^*=\Sigma$.

Throughout, $L=\log d$, logarithms are natural, $x_+=\max(x,0)$,
and $*$ denotes conjugate transpose. Write $I_m$ for the $m\times m$
identity matrix and $\ind_E$ for the indicator of an event or set $E$.
We use the Euclidean norm on
vectors, the operator norm $\|\cdot\|_\op$ on matrices, and the
Hilbert--Schmidt norm $\|A\|_{\rm HS}=(\tr A^*A)^{1/2}$.
For a random variable $Y$, $\|Y\|_q=(\E|Y|^q)^{1/q}$.
The functions $\Gamma$ and $\zeta$ are the gamma and Riemann zeta
functions, and $\gamma$ is Euler's constant.
Positive constants $c,C$ may change between occurrences;
subscripts indicate their permitted parameter dependence.
Our spectral assumption is
\begin{equation}\label{eq:bounded}
 v_d>0,\qquad
 \max_k|\lambda_{k,d}-\bar\lambda_d|\le K\sqrt{v_d},
\end{equation}
where $K$ is independent of $d$.

For the Erd\H{o}s--R\'enyi graph $G(n,t/n)$, in which edges are
independently present with probability $t/n$, the limiting expected
number of components per vertex is (see \citet[Chapter~5]{JLR})
\begin{equation}\label{eq:ct-definition}
 c(t)=\sum_{k\ge1}\frac{k^{k-2}}{k!}t^{k-1}e^{-kt},\qquad t>0.
\end{equation}
The elementary properties of $c$ are proved in
Lemma~\ref{lem:exact-subcritical}.

The classical independent-model constant is
\begin{equation}\label{eq:AB-curve}
 \mathcal C(s)=s\int_0^\infty c(t)t^{s-1}\,dt
 =s\sum_{k\ge1}\frac{\Gamma(k+s-1)}{k!k^{s+1}},\qquad s>0.
\end{equation}
It extends analytically to zero with $\mathcal C(0)=1$.
In the classical independent model with distribution function
proportional to $x^D$ near zero, the parameter is $s=1/D$.
This dimension parameter $D$ is distinct from the matrix size $d$.
For the Haar model, put
\[
 Q_e=\exp(-|H_e|^2/\sigma_d^2),\qquad
 \mathcal C_d(s)=d^{s-1}\E\sum_{e\in\mathcal T_d}Q_e^s,
 \quad s>0,
\]
where $\mathcal T_d$ is the MST edge set selected by the algorithm of
\citet{Kruskal}, which considers edges in increasing order and adds an
edge precisely when it joins different components. Ties are resolved
by a fixed rule. This tree minimizes every
strictly increasing transform of the $Q_e$, including $Q_e^s$ and
$|H_e|^{-1}$. In the iid model the $Q_e$ are independent uniform
variables. Unmarked notation refers to the Haar model; iid analogues
carry an iid superscript when both ensembles appear together.

\begin{theorem}[Expected MST weights under Haar conjugation]\label{thm:power}
Under \eqref{eq:bounded}, for every fixed $S>0$, uniformly for
$0<s\le S$,
\begin{equation}\label{eq:power-main}
 |\mathcal C_d(s)-\mathcal C(s)|
 \le C_{K,S}d^{-1/40}(\log d)^{2s}.
\end{equation}
In particular, at $s=1$,
\begin{equation}\label{eq:frieze-main}
 \mathcal C_d(1)=\E\min_T\sum_{e\in T}e^{-|H_e|^2/\sigma_d^2}
 =\zeta(3)+O_K\bigl(d^{-1/40}(\log d)^2\bigr).
\end{equation}
Both estimates hold for the iid model. Thus the Haar expectation converges to $\mathcal C(s)$ for every fixed $s>0$, including
all dimension parameters $D\ge1$ through $s=1/D$.
\end{theorem}

\paragraph*{Proof outline}
Lemmas~\ref{lem:frame-net}--\ref{lem:external-growing} yield the
component estimate in Theorem~\ref{thm:components}. The conditional
row law in Lemmas~\ref{lem:row}--\ref{lem:radius} controls very small
edge probabilities, while Lemma~\ref{lem:connectivity} gives
connectivity above a polylogarithmic threshold. Integrating these
three estimates in \eqref{eq:power-kruskal} proves the theorem in
Section~\ref{sec:power-laplace}.

The contribution is the uniform transfer to Haar-conjugated matrices.
For independent uniform weights at $s=1$, much sharper finite-size
asymptotics are available from \citet{CooperEtAl}.

Extending $\mathcal C_d$ continuously to $s=0$ gives
the value at zero
\begin{equation}\label{eq:curve-zero}
 \mathcal C_d(0)=\frac{d-1}{d}.
\end{equation}
It also gives an exact Laplace representation of the normalized
mean MST weight for reciprocal edge weights:
\begin{equation}\label{eq:laplace-main}
 F_d=\frac{\sigma_d\sqrt L}{d}\E M_d
 =\sqrt{\frac L\pi}\int_0^\infty
       s^{-1/2}e^{-Ls}\mathcal C_d(s)\,ds.
\end{equation}
The identity holds in both ensembles and is proved in
Section~\ref{sec:power-laplace}. Uniformity of
\eqref{eq:power-main} near zero determines the limiting profile below.

\subsection{Reciprocal edge weights}
For $L>0$, put
\begin{equation}\label{eq:Phi}
 \Phi(L)=1+\frac1{2L}\int_0^{e^{L/2}}
 [c(t)-\ind_{\{t<1\}}]
 \left(1-\frac{\log t}{L}\right)^{-3/2}\frac{dt}{t}.
\end{equation}
This is a convergent, one-dimensional integral. The upper cutoff is
$\sqrt d$ when $L=\log d$; any fixed power $d^\beta$, $0<\beta<1$,
changes the profile by $O(e^{-c d^{\min(\beta,1/2)}})$ for some
$c>0$ and sufficiently large $d$.

\begin{theorem}[Integral law and mean universality]\label{thm:universality}
Under \eqref{eq:bounded},
\begin{equation}\label{eq:integral-main}
 F_d=\frac{\sigma_d\sqrt L}{d}\E M_d
 =\Phi(L)+O_K(d^{-1/40}).
\end{equation}
The same conclusion holds for $M_d^{\mathrm{iid}}$. Consequently,
\begin{equation}\label{eq:universality}
 \left|\frac{\sigma_d\sqrt L}{d}
       (\E M_d-\E M_d^{\mathrm{iid}})\right|
 =O_K(d^{-1/40}).
\end{equation}
\end{theorem}

\paragraph*{Proof outline}
Proposition~\ref{prop:laplace} gives the integral identity and bounds
its part with $s\ge1/2$, using Lemma~\ref{lem:connectivity} and the
inverse moments in Lemma~\ref{lem:radius}. Theorem~\ref{thm:power}
then controls the remaining integral uniformly as $s\downarrow0$.
The proof follows Proposition~\ref{prop:laplace}.

The uniform estimate in Theorem~\ref{thm:power} gives the approximation
in Theorem~\ref{thm:universality} through the exact transform
\eqref{eq:laplace-main}. The profile $\Phi$ keeps the factor
$(1-\log t/L)^{-3/2}$ in \eqref{eq:Phi} without expanding it.
The exponent $1/40$ is limited by the component-size cutoff
$d^{1/40}$ and the error terms in Theorem~\ref{thm:components};
we do not optimize it here.

\subsection{Asymptotic coefficients and a law of large numbers}
Define
\begin{equation}\label{eq:moments-definition}
 \mu_r=\int_0^\infty[c(t)-\ind_{\{t<1\}}](\log t)^r\,\frac{dt}{t},
 \qquad r=0,1,2,\ldots.
\end{equation}
\begin{corollary}[Divergent asymptotic expansion]\label{thm:expansion}
Under \eqref{eq:bounded}, for every fixed integer $N\ge1$,
\begin{equation}\label{eq:expansion}
 F_d=1+\sum_{r=0}^{N-1}\frac{(3/2)_r\mu_r}{2r!L^{r+1}}
        +O_{K,N}(L^{-N-1}).
\end{equation}
The same expansion holds for the iid model and for $\Phi(L)$. Writing
$a_{r+1}=(3/2)_r\mu_r/(2r!)$, one has
\begin{equation}\label{eq:second-order}
 a_1=1-\frac\gamma2-\frac{\pi^2}{12}\approx-0.1111,
 \qquad a_2\approx0.4627,
\end{equation}
and
\begin{equation}\label{eq:factorial-main}
 a_{r+1}\sim\frac{(-1)^{r+1}}4(3/2)_r\qquad(r\to\infty).
\end{equation}
Thus the formal power series in $1/L$ has radius of convergence zero.
Here $(x)_r=x(x+1)\cdots(x+r-1)$, with $(x)_0=1$, is the rising factorial.
\end{corollary}

\paragraph*{Proof outline}
Proposition~\ref{prop:coefficient-evaluation} evaluates the logarithmic
moments of $c$, and Proposition~\ref{prop:coefficient-divergence}
proves the expansion and its factorial coefficient growth.
Combining that expansion with Theorem~\ref{thm:universality} proves
the corollary; the proof is given in Section~\ref{sec:coefficients}.

Each fixed truncation has the remainder in \eqref{eq:expansion}.
For numerical evaluation, we use the integral \eqref{eq:Phi}
instead of truncating its divergent expansion.
Appendix~\ref{sec:numerics} gives a deterministic quadrature check.

The first-order estimate also gives convergence in $L^1$.
\begin{theorem}[Law of large numbers for reciprocal weights]\label{thm:main}
Under \eqref{eq:bounded},
\begin{equation}\label{eq:main}
 \frac{\sqrt{v_d\log d}}{d^{3/2}}M_d\longrightarrow1
 \quad\text{in }L^1.
\end{equation}
\end{theorem}

\paragraph*{Proof outline}
Lemmas~\ref{lem:row}--\ref{lem:max} give the recursive-tree upper
bound in Section~\ref{sec:recursive}. The tail probability for a single
edge, given by Lemma~\ref{lem:row}, together with radius concentration from
Lemma~\ref{lem:radius}, gives a matching lower bound by counting
edges below a suitable weight. Section~\ref{sec:lln-proof} completes
the $L^1$ argument.

\subsection{High-dimensional Euclidean minimum spanning trees}
The spatial dimension in the following statement is denoted by $D$.
The symbol $\omega_D$ denotes unit-ball volume; $v_d$ continues to
denote spectral variance.

\begin{proposition}[Euclidean correction from Penrose's limit]\label{cor:geometric-endpoint}
Let $D\ge2$ be an integer, let
$\omega_D=\pi^{D/2}/\Gamma(1+D/2)$, and let $\beta_D^{\rm geom}$ be the
Euclidean MST constant defined by
\[
 \beta_D^{\rm geom}
 =\lim_{n\to\infty}n^{-(D-1)/D}\E L_n^{(D)},
\]
where $L_n^{(D)}$ is the Euclidean MST length of $n$ independent
uniform points in $[0,1]^D$. Existence of this constant follows from
\citet{SteeleEuclidean}. Then
\begin{equation}\label{eq:penrose-degree-correction}
 \omega_D^{1/D}\beta_D^{\rm geom}
 =\mathcal C(1/D)+o(D^{-1})
 =1+\frac{2-\gamma-\zeta(2)}{D}+o(D^{-1}),
 \qquad D\to\infty.
\end{equation}
The limit in $n$ in the definition of $\beta_D^{\rm geom}$ is taken
at fixed $D$ before the limit $D\to\infty$.
\end{proposition}

\paragraph*{Proof outline}
We combine the component integral of \citet{AvramBertsimas}
with the component limit of \citet{Penrose}. The claim reduces
to convergence of an integral. A bound using the degree
of the origin supplies an integrable bound independent of $D$.
This gives a short proof of the logarithmic-moment consequence of
\citet[Theorem~2 and Lemmas~3--4]{PenroseHighDim} discussed in
Section~\ref{sec:geometric}.

\subsection{A common limiting edge-weight distribution}
Let $\tau$ have survival function $\Prob(\tau>t)=c(t)$, $t>0$.
This is the limiting distribution of $n$ times a uniformly selected
edge weight from the MST with iid uniform weights on $n$ vertices,
identified by \citet{Aldous1990}.
The standard integral formula for moments gives
\[
 \mathcal C(s)=\E\tau^s,\qquad
 \Phi(L)=\E\left(1-\frac{\log(\tau\wedge e^{L/2})}{L}\right)^{-1/2},
 \qquad s,L>0.
\]
For the second identity, integrate the derivative of
$(1-\log t/L)^{-1/2}$ against $c(t)$ on $(0,e^{L/2})$; its value
tends to zero at $t=0$. Moreover,
\[
 a_m=\frac{(1/2)_m}{m!}\E(\log\tau)^m,\qquad m\ge1.
\]
Indeed, integration by parts gives
$\E(\log\tau)^m=m\mu_{m-1}$.
Thus the constants for weights $Q_e^s$, the reciprocal-weight
approximation, and the Euclidean correction all use moments of
the same distribution. Its logarithm has an exponential left tail:
Lemma~\ref{lem:exact-subcritical} gives
$\Prob(\log\tau\le-x)=e^{-x}/2$ for $x\ge0$.
This tail explains the factorial growth of the expansion coefficients.

\subsection{Organization of the proof}
The proof starts from the exact conditional Dirichlet law of one
off-diagonal row. It gives the first-order estimate through an explicit recursive
tree. For the weights $Q_e^s$, we control components with up to $d^{1/40}$ vertices uniformly for
$d^{-1/40}\le t\le d^{1/40}$. The proof represents their internal edges
using independent Gaussian vectors and represents their external edges by a
Gaussian polar decomposition. The resulting estimates are relative to
the probabilities of the component events, which can be exponentially
small in the component size.

Gaussian approximation of sphere coordinates is treated by
\citet{DiaconisFreedman}. For Haar orthogonal and unitary matrices,
\citet{JiangMaxima} studies simultaneous entrywise approximation and
the largest entries, while \citet{Jiang} and \citet{JiangUnitary}
treat growing blocks in the orthogonal and unitary cases, respectively.
For unitary conjugates of deterministic spectra,
\citet[Theorem~4.1]{OlshanskiVershik} characterize the limits of
fixed corners; under \eqref{eq:bounded}, the limit is Gaussian.
The estimates below concern joint edge-presence and edge-absence
events for growing vertex sets, at probabilities small enough that
absolute approximation errors do not suffice.
The global tail estimate comes from
Haar concentration applied to the distance of the matrix entries across a graph cut from
the set of blocks with all entries below threshold. It yields
connectivity at mean-degree parameter $A_K(\log d)^2$, as shown
in Section~\ref{sec:power-laplace}.

\section{The conditional row distribution and inverse moments}
Write $u_j$ for column $j$ of $U_d$ and put
\[
 R_j^2=u_j^*\Lambda_d^2u_j-(u_j^*\Lambda_du_j)^2.
\]
Thus $R_j$ is the Euclidean norm of the off-diagonal entries in
row $j$. This quantity is unchanged by translating the spectrum.

\begin{lemma}\label{lem:row}
Conditional on $u_j$, the vector
\[
 \left(\frac{|H_{ij}|^2}{R_j^2}\right)_{i\ne j}
 \quad\text{has law }\operatorname{Dirichlet}(1,\ldots,1).
\]
Its conditional distribution is independent of $u_j$; consequently the
normalized vector is independent of $R_j$.
\end{lemma}
\begin{proof}
Conditional on $u_j$, the other columns form a Haar orthonormal basis of
$u_j^\perp$. The off-diagonal entries are the coordinates, in that basis,
of the fixed vector $(I-u_ju_j^*)\Lambda_du_j$, whose norm is $R_j$.
Its normalized coordinates are uniform on the complex unit sphere in
$\mathbb C^{d-1}$. Their squared moduli are normalized independent
exponentials, giving the assertion. For a nonconstant spectrum, $R_j>0$
almost surely.
\end{proof}

\begin{lemma}\label{lem:radius}
Under \eqref{eq:bounded}, for $0<\delta<1$,
\begin{equation}\label{eq:concentration}
 \Prob\left(\left|\frac{R_j^2}{v_d}-1\right|>\delta\right)
 \le C_K\exp(-c_Kd\delta^2).
\end{equation}
Moreover, as $d\to\infty$,
\begin{equation}\label{eq:inv}
 \E R_j^{-1}\sim v_d^{-1/2}.
\end{equation}
For every nonconstant spectrum, without assumption \eqref{eq:bounded},
\begin{equation}\label{eq:inverse-radius-general}
 \E(R_j^2)^{-q}
 \le\frac{\Gamma(d)}{\Gamma(d-2q)}
       \frac{\Gamma(1-q)^2}{(d^2v_d)^q},\qquad 0<q<1.
\end{equation}
\end{lemma}
\begin{proof}
Set $b_k=(\lambda_{k,d}-\bar\lambda_d)/\sqrt{v_d}$, so
$d^{-1}\sum b_k=0$, $d^{-1}\sum b_k^2=1$, and $|b_k|\le K$.
The column weights $P_k=|u_{kj}|^2$ have the representation
$P_k=E_k/S_d$, where the $E_k$ are independent mean-one exponentials and
$S_d=\sum_kE_k$. Therefore
\begin{equation}\label{eq:rform}
 \frac{R_j^2}{v_d}=\sum_k b_k^2P_k-\left(\sum_k b_kP_k\right)^2.
\end{equation}
For bounded real coefficients $c_k$, the exponential moment identity
\[
 \log\E\exp\left(\theta\sum_kc_k(E_k-1)\right)
 =\sum_k[-\theta c_k-\log(1-\theta c_k)]
 \le C\theta^2\sum_kc_k^2
\]
holds when $|\theta|\max|c_k|\le1/2$. Chernoff bounds applied to
$c_k=1,b_k,b_k^2$, control $S_d/d$, $\sum b_kE_k/d$, and
$\sum b_k^2E_k/d$ with errors $O(\delta)$ and failure probability
$C_K e^{-c_Kd\delta^2}$. Equation \eqref{eq:rform} proves
\eqref{eq:concentration}, after adjusting constants.

To prove the negative-moment estimate, let
$X_1,\ldots,X_d$ be independent mean-one exponentials and $S=\sum_iX_i$.
Put
\[
 \Delta(X)=\sum_{i<j}(\lambda_{i,d}-\lambda_{j,d})^2X_iX_j,
 \qquad W=\sum_{i<j}(\lambda_{i,d}-\lambda_{j,d})^2=d^2v_d.
\]
Then $R_j^2=\Delta(X)/S^2$ in distribution, and $S$ is independent
of the normalized Dirichlet weights determining $R_j$. For $0<q<1$,
convexity of $x\mapsto x^{-q}$ and normalized weights
$(\lambda_{i,d}-\lambda_{j,d})^2/W$ give
\[
 \E \Delta(X)^{-q}\le W^{-q}\Gamma(1-q)^2.
\]
On the other hand,
$\E \Delta(X)^{-q}=\E S^{-2q}\,\E(R_j^2)^{-q}$, so
\eqref{eq:inverse-radius-general} follows.
This bound holds for every nonconstant spectrum. The special case
$q=1/2$ was recorded by \citet[Equation~(17)]{CMcontrol}:
\[
 \E R_j^{-1}\le\frac{\pi(d-1)}{d\sqrt{v_d}}.
\]
Taking $q=3/4$ in \eqref{eq:inverse-radius-general} bounds the
$3/2$ moment of $\sqrt{v_d}/R_j$ uniformly in $d$. Combined with
\eqref{eq:concentration}, this gives uniform integrability and
proves \eqref{eq:inv}.
\end{proof}

\section{The recursive tree and its expected weight}\label{sec:recursive}
For each $j=2,\ldots,d$, choose a parent
\[
 p(j)\in\mathop{\rm arg\,max}_{1\le i<j}|H_{ij}|.
\]
Break ties deterministically. These edges form a spanning tree: following
parent pointers strictly decreases the label and ends at $1$. Thus
\begin{equation}\label{eq:tree}
 M_d\le M_d^{\rm rec}:=\sum_{j=2}^d\left(\max_{i<j}|H_{ij}|\right)^{-1}.
\end{equation}

Let $n=d-1$, let $\xi=(\xi_1,\ldots,\xi_n)$ have the Dirichlet law in
Lemma \ref{lem:row}, and set $Z_m=\max(E_1,\ldots,E_m)$ for
independent mean-one exponentials $E_1,\ldots,E_n$.
With $S_n=\sum_{i=1}^nE_i$, the gamma--Dirichlet representation
$E_i=S_n\xi_i$ has $S_n$ independent of $\xi$. Consequently,
\begin{equation}\label{eq:gamma}
 \E\left(\max_{i\le m}\xi_i\right)^{-1/2}
 =\frac{\Gamma(n)}{\Gamma(n-1/2)}\E Z_m^{-1/2},
 \qquad 1\le m\le n.
\end{equation}
The identity also holds for $m=1$.

\begin{lemma}\label{lem:max}
Writing $h_m=\E Z_m^{-1/2}$, one has
\[
 h_m\sim(\log m)^{-1/2},\qquad
 \sum_{m=1}^n h_m\sim\frac n{\sqrt{\log n}}.
\]
\end{lemma}
\begin{proof}
The distribution function of $Z_m$ is $(1-e^{-x})^m$. For fixed
$q>0$ and $m>q$,
\[
 \E Z_m^{-q}=q\int_0^\infty x^{-q-1}(1-e^{-x})^m\,dx.
\]
Fix $0<\alpha<1$ and split at $1$ and $\alpha\log m$ (for large $m$).
The first integral is at most $q/(m-q)$, using $1-e^{-x}\le x$.
The second is at most $e^{-m^{1-\alpha}}$, using
$(1-e^{-x})^m\le e^{-me^{-x}}$. The last is at most
$(\alpha\log m)^{-q}$. For the reverse bound,
$\Prob(Z_m\le(1+\varepsilon)\log m)\to1$ by a union bound.
Let $\alpha\uparrow1$ and $\varepsilon\downarrow0$ to get
$\E Z_m^{-q}\sim(\log m)^{-q}$.

For the sum, $h_m$ decreases, so its lower bound is $nh_n$.
The first $n^\alpha$ terms are at most $n^\alpha h_1$; the rest are
at most $n(1+o(1))/\sqrt{\alpha\log n}$. Let $\alpha\uparrow1$.
\end{proof}

Taking expectations in \eqref{eq:tree}, using the identical marginal
laws of the radii and Lemma \ref{lem:row}, gives
\begin{equation}\label{eq:ET}
 \E M_d^{\rm rec}=\E R_1^{-1}
 \frac{\Gamma(d-1)}{\Gamma(d-3/2)}
 \sum_{m=1}^{d-1}h_m.
\end{equation}
Equations \eqref{eq:inv} and \eqref{eq:ET}, Lemma~\ref{lem:max}, and the gamma ratio asymptotic
therefore show
\begin{equation}\label{eq:upper}
 \E M_d\le\E M_d^{\rm rec}\sim\frac{d^{3/2}}{\sqrt{v_d\log d}}.
\end{equation}

There is also a universal upper bound without \eqref{eq:bounded}.
For every nonconstant deterministic spectrum, the inverse-radius
estimate \eqref{eq:inverse-radius-general} at $q=1/2$ gives
$\E R_1^{-1}\le\pi(d-1)/(d\sqrt{v_d})$. Substitution into
\eqref{eq:ET} yields
\[
 \E M_d\le(\pi+o(1))\frac{d^{3/2}}{\sqrt{v_d\log d}}.
\]
The $o(1)$ here is independent of the spectrum. This proves the
logarithmic refinement suggested by \citet[Section~1.4]{CMcontrol}
of their path-tree bound in Theorem~11, which also applies to every
nonconstant spectrum.

\section{Matching lower bound and \texorpdfstring{$L^1$}{L1} convergence}\label{sec:lln-proof}
\begin{proof}[Proof of Theorem~\ref{thm:main}]
Fix $\varepsilon>0$ and define
\[
 u_{d,\varepsilon}^2=(1+\varepsilon)\frac{v_d\log d}{d},\qquad
 N_d=\sum_{i<j}\ind_{\{|H_{ij}|>u_{d,\varepsilon}\}}.
\]
Conditional on $u_1$, Lemma \ref{lem:row} gives the exact beta tail
\begin{equation}\label{eq:beta}
 \Prob(|H_{12}|>u_{d,\varepsilon}\mid u_1)
 =\left(1-\frac{u_{d,\varepsilon}^2}{R_1^2}\right)_+^{d-2}.
\end{equation}
Choose a fixed $0<\delta<\varepsilon$. On
$R_1^2\le(1+\delta)v_d$, \eqref{eq:beta} is at most
\[
 \exp\left[-\frac{d-2}{d}\frac{1+\varepsilon}{1+\delta}\log d\right]
 =d^{-(1+\varepsilon)/(1+\delta)+o(1)}=o(d^{-1}).
\]
The complementary event has probability $O(e^{-c d})$ by
Lemma \ref{lem:radius}. It follows that
\begin{equation}\label{eq:N}
 \E N_d=\binom d2\Prob(|H_{12}|>u_{d,\varepsilon})=o(d).
\end{equation}
Every spanning tree contains $d-1$ edges, and at most $N_d$ of them
can have weight less than $u_{d,\varepsilon}^{-1}$. Therefore
\begin{equation}\label{eq:lower}
 M_d\ge\frac{d-1-N_d}{u_{d,\varepsilon}}.
\end{equation}
The inequality remains true when its right-hand side is negative.
Taking expectations in \eqref{eq:lower}, using \eqref{eq:N}, and
then letting $\varepsilon\downarrow0$ gives the matching lower
bound for \eqref{eq:upper}.

Let $\bar M_d=\sqrt{v_d\log d}\,M_d/d^{3/2}$.
We have proved $\E \bar M_d\to1$. Also $N_d/d\to0$ in probability by
Markov's inequality, so \eqref{eq:lower} implies
$\Prob(\bar M_d<1-\eta)\to0$ for every $\eta>0$. Since
$0\le(1-\bar M_d)_+\le1$, its expectation tends to zero. Hence
\[
 \E|\bar M_d-1|=\E \bar M_d-1+2\E(1-\bar M_d)_+\longrightarrow0.
\]
\end{proof}

\section{Quantitative universality of expected component counts}
\label{sec:component-universality}

Let $\mathcal G_d(t)$ be the graph
whose edges satisfy
\[
 |H_{ij}|^2>\sigma_d^2\log(d/t),\qquad 0<t<d.
\]
Write $\kappa_d(t)$ for its number of connected components. Throughout
this section the standardized diagonal matrix
\[
 A_d=\frac{\Lambda_d-\bar\lambda_dI}{\sqrt{v_d}}
 \quad\text{satisfies}\quad
 \tr A_d=0,\quad d^{-1}\tr A_d^2=1,\quad\|A_d\|_\op\le K,
\]
where $K$ does not depend on $d$. The component estimate below is uniform on a growing interval of edge
probabilities and for components of growing size.

\begin{theorem}\label{thm:components}
Set $\alpha=1/40$ and write $c_k(t)$ for the $k$th term of
the series \eqref{eq:ct-definition}, so that
\begin{equation}\label{eq:component-series}
 c_k(t)=\frac{k^{k-2}}{k!}t^{k-1}e^{-kt},\qquad
 c(t)=\sum_{k\ge1}c_k(t).
\end{equation}
Then
\begin{equation}\label{eq:component-rate}
 \sup_{d^{-\alpha}\le t\le d^\alpha}
 \left|\frac{\E\kappa_d(t)}d-c(t)\right|
 =O_K(d^{-\alpha}).
\end{equation}
The same estimate holds for the component count of the
Erd\H{o}s--R\'enyi graph $G(d,t/d)$, whose edges are independently
present with probability $t/d$.
Consequently the normalized expected component counts in the two
models differ by $O_K(d^{-\alpha})$, uniformly on this interval.
\end{theorem}

\citet[Proposition~2]{CM} establish a qualitative Gaussian approximation
for fixed finite collections of off-diagonal entries under a
fourth-moment Lyapunov condition.
After centering the spectrum, \eqref{eq:bounded} implies that condition,
since $\sum_k(\lambda_{k,d}-\bar\lambda_d)^4/(d^2v_d^2)\le K^2/d$.
Here the component size and thresholds grow with $d$, and the bounds
are relative to the rare-event probabilities.
We prove them using the following Gaussian representations. Put
\[
 L=\log d,\qquad J=\lfloor d^\alpha\rfloor,\qquad
 \delta=d^{-1/8},\qquad p=t/d,\qquad \ell=\log(d/t).
\]
On the interval in \eqref{eq:component-rate},
$(1-\alpha)L\le\ell\le(1+\alpha)L$. Constants below may depend on
$K$ but are uniform in $k\le J$ and in $t$ in this interval.

\begin{lemma}[Concentration of a growing Haar frame]\label{lem:frame-net}
Let $V=(u_1,\ldots,u_k)$ consist of $k\le J$ Haar columns. Outside
an event of probability at most
\begin{equation}\label{eq:frame-failure}
 C\exp(Ck-cd\delta^2),
\end{equation}
one has
\begin{equation}\label{eq:frame-good}
 \|V^*A_dV\|_\op\le C\delta,
 \qquad
 \|V^*A_d^2V-I_k\|_\op\le C\delta.
\end{equation}
We call this event $\mathcal F_k$.
\end{lemma}
\begin{proof}
For any fixed unit vector $x\in\mathbb C^k$, $Vx$ is uniform on
the complex unit sphere in $\mathbb C^d$. Its squared coordinate
moduli have the representation $E_i/\sum_jE_j$, with independent
mean-one exponential variables. If $|b_i|\le B$, the identity
\[
 \log\E\exp\left(\theta\sum_i b_i(E_i-1)\right)
 =\sum_i[-\theta b_i-\log(1-\theta b_i)]
 \le C_Bd\theta^2
\]
holds for $|\theta|\le(2B)^{-1}$. The same exponential Markov inequality and the
same estimate for $\sum E_i$ therefore imply, for $0<s<1$,
\[
 \Prob\left(\left|(Vx)^*\Delta(Vx)-d^{-1}\tr \Delta\right|>s\right)
 \le C_Be^{-c_Bds^2}
\]
for real diagonal $\Delta$ with $\|\Delta\|_\op\le B$.
Apply this to $\Delta=A_d$ and $\Delta=A_d^2$. A $1/4$-net of the unit
sphere in $\mathbb C^k$ has at most $9^{2k}$ points. For every
Hermitian $k\times k$ matrix $Q$, its operator norm is at most
twice the maximum of $|x^*Qx|$ on this net: approximate a unit
maximizer by a net point and bound the two resulting error terms
by $\|Q\|_\op/2$. Taking a union bound proves
\eqref{eq:frame-failure}--\eqref{eq:frame-good}.
\end{proof}

\begin{lemma}[Relative moderate deviations for internal edges]
\label{lem:internal-growing}
Let $F$ be any set of $r\le k$ distinct edges among $k\le J$
specified vertices, and write $E(G)$ for the edge set of a graph $G$.
Uniformly over their labels and their edge pattern,
\begin{equation}\label{eq:internal-growing}
 \Prob(F\subseteq E(\mathcal G_d(t)))
 =p^r\left[1+O_K(d^{-1/16}L)\right].
\end{equation}
\end{lemma}
\begin{proof}
By exchangeability take the vertices to be $1,\ldots,k$.
Generate their Haar columns sequentially from independent
$g_j\sim\CN(0,I_d)$ by
\[
 u_j=\frac{P_{j-1}g_j}{\nu_j},\qquad
 P_{j-1}=I-V_{j-1}V_{j-1}^*,\qquad \nu_j=\|P_{j-1}g_j\|,
\]
where $V_{j-1}=(u_1,\ldots,u_{j-1})$. Let
\[
 B_j=P_{j-1}A_dV_{j-1},\qquad
 \Sigma_j=B_j^*B_j
 =V_{j-1}^*A_d^2V_{j-1}-(V_{j-1}^*A_dV_{j-1})^2.
\]
Conditional on the past, meaning all Gaussian and auxiliary vectors
generated at earlier steps, $B_j^*g_j$ is complex Gaussian with
covariance $\Sigma_j$. Define a standard complex Gaussian vector by
\[
 z_j=\Sigma_j^{+1/2}B_j^*g_j+\Pi_{\ker \Sigma_j}\xi_j,
\]
where $\Sigma_j^{+1/2}$ denotes the positive square root of the
Moore--Penrose inverse, $\Pi_{\ker \Sigma_j}$ is the orthogonal projection
onto its kernel, and $\xi_j\sim\CN(0,I_{j-1})$ is an independent
auxiliary vector. Then, even if $\Sigma_j$ is singular,
\[
 z_j\mid\text{past}\sim\CN(0,I_{j-1}),\qquad
 B_j^*g_j=\Sigma_j^{1/2}z_j.
\]
Thus all coordinates of all the $z_j$ are mutually independent
standard complex Gaussian variables.

On $\mathcal F_k$, every principal subframe obeys
\eqref{eq:frame-good}, and hence
$\|\Sigma_j-I_{j-1}\|_\op\le C\delta$ for $j\le k$.
Conditional on the past, $\nu_j^2$ is gamma with shape $d-j+1$
and rate one.
The exponential moment bound used above gives
\[
 \Prob\left(\left|\frac{\nu_j^2}{d-j+1}-1\right|>\delta
       \;\middle|\;\text{past}\right)
 \le Ce^{-cd\delta^2}.
\]
Except on an event of probability $C\exp(Ck-cd\delta^2)$,
all these norm bounds and $\mathcal F_k$ hold simultaneously.
The vector of normalized internal entries in column $j$ is then
\[
 \left(\frac{H_{ij}}{\sigma_d}\right)_{i<j}
 =\Theta_j z_j,\qquad
 \Theta_j=\frac{\sqrt{v_d}}{\sigma_d\nu_j}\Sigma_j^{1/2},\qquad
 \|\Theta_j-I_{j-1}\|_\op\le C(\delta+k/d)\le C\delta.
\]
Here $\sigma_d^2/v_d=d/(d^2-1)$ accounts for the deterministic
normalizing factor.

Impose also $\max_{i,j}|z_{j,i}|^2\le10kL$. Its failure
probability is at most $k^2d^{-10k}$. On this event
$\|z_j\|\le\sqrt{10}\,k\sqrt L$, so every internal coordinate
satisfies
\begin{equation}\label{eq:internal-coordinate-error}
 \left|\frac{|H_{ij}|^2}{\sigma_d^2}-|z_{j,i}|^2\right|
 \le C\delta k^{3/2}L=:\eta_{\rm in}.
\end{equation}
The quadratic error in $\Theta_j-I$ is absorbed here because
$\delta\sqrt k\to0$ uniformly for $k\le J$.

Sandwich the $r$ selected exceedances between the independent
Gaussian exceedance events at levels $\ell+\eta_{\rm in}$ and
$\ell-\eta_{\rm in}$. Their probabilities are
$p^r e^{\mp r\eta_{\rm in}}$, and
\[
 r\eta_{\rm in}\le C\delta k^{5/2}L
 \le Cd^{-1/16}L=o(1).
\]
The discarded probabilities remain negligible after division by
$p^r$. Indeed, $p^r\ge d^{-(1+\alpha)k}$, whereas
\[
 \frac{k^2d^{-10k}}{p^r}\le k^2d^{-(9-\alpha)k},\qquad
 \frac{C\exp(Ck-cd^{3/4})}{p^r}
 \le C\exp\{Ck+(1+\alpha)kL-cd^{3/4}\}.
\]
Both are smaller than the displayed error in
\eqref{eq:internal-growing}. The case $r=0$ is immediate.
\end{proof}

\begin{lemma}[Uniform conditional absence of external edges]
\label{lem:external-growing}
Conditional on any frame $V=(u_1,\ldots,u_k)$ in $\mathcal F_k$,
the probability that no edge joins $\{1,\ldots,k\}$ to its
complement is
\begin{equation}\label{eq:external-growing}
 (1-p)^{k(d-k)}\left[1+O_K(d^{-1/16}L)\right],
\end{equation}
uniformly for $k\le J$ and $d^{-\alpha}\le t\le d^\alpha$.
\end{lemma}
\begin{proof}
Put $n=d-k$ and let $B=(H_{ij})_{i\le k<j}$ be the external
$k\times n$ block. Conditional on $V$, the remaining columns are
a Haar orthonormal basis of $\operatorname{range}(V)^\perp$.
The residual Gram matrix is
\[
 G=v_d\{V^*A_d^2V-(V^*A_dV)^2\}.
\]
On $\mathcal F_k$ it is positive definite and
$\|G/(n\sigma_d^2)-I_k\|_\op\le C\delta$.
If $Z$ is a $k\times n$ matrix of independent $\CN(0,1)$
variables, independent of $V$, the exact conditional law is
\[
 \frac{B}{\sigma_d}\ \stackrel{\rm law}= C_ZZ,
 \qquad
 C_Z=\left(\frac{G}{n\sigma_d^2}\right)^{1/2}
       \left(\frac{ZZ^*}{n}\right)^{-1/2}.
\]
This is the Gaussian polar representation of a uniform
row-orthonormal matrix, followed by multiplication by $G^{1/2}$.

For a fixed unit vector $x\in\mathbb C^k$, the squared coordinates
of $x^*Z$ are independent mean-one exponentials. Their moment
bound and the same $1/4$-net argument as in
Lemma \ref{lem:frame-net} give
\[
 \Prob\left(\left\|ZZ^*/n-I_k\right\|_\op>C\delta\right)
 \le C\exp(Ck-cn\delta^2)=:\epsilon_d.
\]
On the complementary event,
$\|C_Z-I_k\|_\op+\|C_Z^{-1}-I_k\|_\op\le C\delta$.
These inequalities give a deterministic sandwich for the event
that all external entries are below threshold. To see this,
if all coordinates of a column $z$ are bounded by $\sqrt u$,
then
\[
 |(C_Zz)_i|\le\sqrt u+C\delta\|z\|
 \le(1+C\delta\sqrt k)\sqrt u.
\]
The reverse implication follows by applying the corresponding
bound to $C_Z^{-1}$. Hence, with
$\eta_{\rm out}=C\delta\sqrt k\,\ell$, the conditional
no-external-edge probability is bounded below and above by
\[
 (1-e^{-(\ell-\eta_{\rm out})})^{kn}-\epsilon_d
 \quad\text{and}\quad
 (1-e^{-(\ell+\eta_{\rm out})})^{kn}+\epsilon_d,
\]
respectively. Both implications hold whenever
$\|ZZ^*/n-I_k\|_\op\le C\delta$.

Since $p\to0$ and $\eta_{\rm out}\to0$ uniformly,
\[
 kn\left|\log(1-pe^{\pm\eta_{\rm out}})-\log(1-p)\right|
 \le Ckt\eta_{\rm out}
 \le C\delta k^{3/2}tL\le Cd^{-1/16}L.
\]
Also $(1-p)^{kn}\ge e^{-2kt}$ for large $d$, and
\[
 \frac{\epsilon_d}{(1-p)^{kn}}
 \le C\exp(Ck-cd^{3/4}+2kt)
 \le C\exp(-c'd^{3/4}),
\]
because $kt\le d^{2\alpha}$. This proves
\eqref{eq:external-growing}.
\end{proof}

\begin{proof}[Proof of Theorem \ref{thm:components}]
First combine the internal and external estimates carefully, since
the internal edge event is rare. Let $A_F$ denote the event that
the $r\le k$ edges of $F$ are present, and let $Q_0=(1-p)^{k(d-k)}$.
On $\mathcal F_k$, Lemma \ref{lem:external-growing} bounds the
conditional external-absence probability by $Q_0(1+O(\eta_d))$,
where $\eta_d=d^{-1/16}L$. Integrating over $A_F$ and using
Lemma \ref{lem:internal-growing} gives
\begin{equation}\label{eq:internal-external-joint}
 \Prob(A_F,\ \text{no external edge})
 =p^rQ_0[1+O_K(\eta_d)].
\end{equation}
Indeed, the error from frames outside $\mathcal F_k$, even divided
by the rare target probability, is at most
\[
 \frac{C\exp(Ck-cd^{3/4})}{p^rQ_0}
 \le C\exp\{Ck+(1+\alpha)kL+2kt-cd^{3/4}\}
 =O(e^{-c'd^{3/4}}).
\]
The resulting error is negligible relative to $p^rQ_0$ uniformly
on the stated range.

Fix a labelled tree $T$ on a specified set of $k$ vertices. Apply
\eqref{eq:internal-external-joint} first to its $k-1$ edges and
then to those edges plus any one of the
$q_k=\binom{k}{2}-(k-1)$ possible additional edges. A union bound
for these additional edges shows that the probability that the
specified set is a component whose induced graph is exactly $T$
equals
\begin{equation}\label{eq:exact-tree-component}
 p^{k-1}(1-p)^{k(d-k)+q_k}
 \left[1+O_K(\eta_d+k^2p)\right].
\end{equation}
Here $(1-p)^{q_k}=1+O(k^2p)$, so incorporating this factor
preserves the relative error. The estimate is uniform over all
labelled trees with $k\le J$.

Let $N_{d,k}^{\rm tr}$ count tree components of size $k$.
There are $k^{k-2}$ labelled trees on a specified vertex set.
Their exact induced-graph events are disjoint, so
\[
 \frac{\E N_{d,k}^{\rm tr}}d
 =\frac1d\binom dk k^{k-2}p^{k-1}
   (1-p)^{k(d-k)+q_k}
   [1+O_K(\eta_d+k^2p)].
\]
The finite-$d$ expression before the bracket equals
\begin{equation}\label{eq:finite-tree-density}
 c_k(t)\left[1+O\left(\frac{k^2}{d}
                  +\frac{k^2t}{d}+\frac{kt^2}{d}\right)\right].
\end{equation}
This follows by expanding
$\prod_{i=0}^{k-1}(1-i/d)$ and
$\{k(d-k)+q_k\}\log(1-t/d)$; their logarithmic errors have the
displayed bounds. Uniformly on the present range all these errors,
as well as $k^2p$, are $O(d^{3\alpha-1})$, hence negligible
relative to $\eta_d$.

Let $N_{d,k}^{\rm cy}$ count connected components of size $k$ with
at least one cycle. Every such graph contains a spanning tree and one
additional edge. Summing \eqref{eq:internal-external-joint} over
these choices, without imposing absence of any other internal
edge, gives
\begin{equation}\label{eq:cyclic-component-bound}
 \frac{\E N_{d,k}^{\rm cy}}d\le C k^2p\,c_k(t),
 \qquad 3\le k\le J.
\end{equation}
Overcounting a cyclic graph is harmless for this upper bound.

The series in \eqref{eq:component-series} has an elementary
uniform summability bound. For $k\ge2$, maximize
$t^{k-1}e^{-kt}$ at $t=(k-1)/k$, and use
$k!\ge c\sqrt k\,(k/e)^k$. This gives
\[
 \sup_{t>0}c_k(t)\le Ck^{-5/2},\qquad
 \sup_{t>0}\sum_{k>J}c_k(t)\le CJ^{-3/2},\qquad
 \sup_{t>0}\sum_{k\ge1}c_k(t)<\infty.
\]
Finally, any graph on $d$ vertices has at most $d/(J+1)$
components larger than $J$. Summing the tree estimates and
\eqref{eq:cyclic-component-bound}, and then using this deterministic
truncation, proves
\[
 \left|\frac{\E\kappa_d(t)}d-c(t)\right|
 \le C_K\left(J^{-1}+\eta_d+J^2p+J^{-3/2}\right)
 =O_K(d^{-\alpha}).
\]
This is uniform in the stated interval. For $G(d,t/d)$ the exact
tree-component formula is the finite-$d$ expression in
\eqref{eq:finite-tree-density}; the cyclic-component bound
\eqref{eq:cyclic-component-bound} and the same truncation apply directly. They prove
the asserted Erd\H{o}s--R\'enyi estimate as well.
\end{proof}

\section{Expected MST weights and their integral representation}\label{sec:power-laplace}
Put
\[
 X_{ij}=|H_{ij}|^2/\sigma_d^2,\qquad Q_{ij}=e^{-X_{ij}},
 \qquad
 \mathcal C_d(s)=d^{s-1}\E\sum_{e\in\mathcal T_d}Q_e^s,
 \quad s>0,
\]
where $\mathcal T_d$ is the MST edge set. Throughout this section,
$R:=R_1$ is the Euclidean norm of the off-diagonal entries in the first row,
as in Lemma~\ref{lem:row}.
All the functions
$Q\mapsto Q^s$ and $Q\mapsto(-\log Q)^{-1/2}$ are strictly
increasing, so Kruskal's algorithm selects the same edge set for these
weights. Write
\[
 f_d(t)=\frac{\E\kappa_d(t)-1}{d},\quad 0<t<d,
 \qquad \mathcal C(s)=s\int_0^\infty c(t)t^{s-1}\,dt.
\]
The component-count identity is
\begin{equation}\label{eq:power-kruskal}
 \mathcal C_d(s)=s\int_0^d f_d(t)t^{s-1}\,dt.
\end{equation}
The function $\mathcal C$ is finite for every $s>0$, with
$\mathcal C(1)=\zeta(3)$. Its regularized Mellin representation
is developed in Proposition~\ref{prop:coefficient-evaluation}.

\begin{lemma}[Connectivity]\label{lem:connectivity}
Under \eqref{eq:bounded}, for every fixed $B>0$ there is
$A=A(K,B)>0$ such that, for all sufficiently large $d$,
\begin{equation}\label{eq:connectivity}
 \Prob\bigl(\mathcal G_d(A(\log d)^2)
                  \text{ is disconnected}\bigr)\le d^{-B}.
\end{equation}
The same estimate holds for the iid threshold graph.
\end{lemma}
\begin{proof}
Write $L=\log d$ and $X_{ij}=|H_{ij}|^2/\sigma_d^2$.
We first record a marginal lower bound. For $1\le\ell\le L$,
the exact row law gives
\[
 \Prob(X_{12}\ge\ell+1\mid R)
 =\left(1-\frac{\sigma_d^2(\ell+1)}{R^2}\right)_+^{d-2}.
\]
Put $\delta=d^{-1/8}$. On the event
$R^2\ge(1-\delta)v_d$, whose probability is
$1-O_K(e^{-c_Kd^{3/4}})$ by Lemma~\ref{lem:radius}, set
$q=\sigma_d^2(\ell+1)/R^2$. Then
$q\le(1+2\delta)(\ell+1)/d$ for large $d$.
The inequality $\log(1-q)\ge-q-2q^2$, valid for $0\le q\le1/2$,
therefore gives
\[
 (1-q)^{d-2}\ge
 \exp\left[-(\ell+1)-O\left(\delta L+L^2/d\right)\right].
\]
Consequently, uniformly for $1\le\ell\le L$,
\begin{equation}\label{eq:cut-marginal}
 \Prob(X_{12}\ge\ell+1)\ge c e^{-\ell}
\end{equation}
for all sufficiently large $d$.

For a nonempty set $S\subset\{1,\ldots,d\}$, $k=|S|\le d/2$,
and $\ell=\log(d/t)$, define
\[
 F_S(U)=\left[\sum_{i\in S,\,j\notin S}
       \left(\frac{|(U^*\Lambda_dU)_{ij}|}{\sigma_d}-\sqrt\ell\right)_+^2
                  \right]^{1/2}.
\]
As a function of the normalized cross block, this is the Euclidean
distance to the product of complex discs of radius $\sqrt\ell$,
so it is $1$-Lipschitz. Centering the spectrum and using
\eqref{eq:bounded} shows that its Lipschitz constant on $U(d)$,
with the Hilbert--Schmidt metric, is at most
\[
 \frac{2\|\Lambda_d-\bar\lambda_d I\|_\op}{\sigma_d}
 \le2K\sqrt{\frac{d^2-1}{d}}\le2K\sqrt d.
\]
The Haar concentration inequality of
\citet[Theorem~5.17]{MeckesBook} states that an $L_F$-Lipschitz
function satisfies
$\Prob(F-\E F\ge r)\le
\exp[-(d-2)r^2/(24L_F^2)]$.
Applying it to $F_S$ and $-F_S$, and integrating the resulting
two-sided tail, gives constants independent of $S,t,d$ such that
\begin{equation}\label{eq:cut-concentration}
 \Prob(F_S\le\E F_S-r)\le e^{-c_Kr^2},
 \qquad \operatorname{Var}(F_S)\le C_K.
\end{equation}

If $X_{ij}\ge\ell+1$, its contribution to $F_S^2$ is at least
$(\sqrt{\ell+1}-\sqrt\ell)^2\ge1/(8L)$.
Exchangeability and \eqref{eq:cut-marginal} imply, for
$1\le\ell\le L$,
\[
 \E F_S^2\ge\frac{c}{L}k(d-k)\frac td
             \ge\frac{ckt}{L}.
\]
Take $t=A L^2$. For fixed $A>0$ and sufficiently large $d$,
$1\le\ell\le L$, and this lower bound dominates the variance
bound uniformly for $k\ge1$. Hence
\[
 (\E F_S)^2=\E F_S^2-\operatorname{Var}(F_S)
       \ge ckt/L.
\]
An empty cut implies $F_S=0$, so \eqref{eq:cut-concentration}
yields
\[
 \Prob(S\text{ has no edge to }S^c)
       \le\exp(-c_Kkt/L)=d^{-c_KAk}.
\]
Every disconnected graph has an empty cut with $k\le d/2$.
Choose $A=A(K,B)$ with $c_KA\ge B+3$ and use
$\binom dk\le d^k$ to obtain
\[
 \Prob(\mathcal G_d(A L^2)\text{ disconnected})
 \le\sum_{k=1}^{\lfloor d/2\rfloor}d^{-(B+2)k}
 \le d^{-B}
\]
for large $d$.
For the iid graph, the probability of a fixed empty cut is
$(1-t/d)^{k(d-k)}\le e^{-kt/2}$ for $k\le d/2$.
The same union bound proves the assertion.
\end{proof}

\begin{lemma}[Properties of $c$]\label{lem:exact-subcritical}
The component density satisfies
\begin{equation}\label{eq:exact-subcritical}
 c(t)=1-\frac t2\quad(0<t\le1),\qquad
 0<c(t)<\frac1{2t}\quad(t>1),\qquad
 c(t)\le e^{1-t}\quad(t\ge1).
\end{equation}
\end{lemma}
\begin{proof}
Let
\[
 \mathsf Y(z)=\sum_{k\ge1}\frac{k^{k-1}}{k!}z^k,
 \qquad
 \mathsf U(z)=\sum_{k\ge1}\frac{k^{k-2}}{k!}z^k.
\]
The rooted-tree generating function (see \citet[Section~II.5]{FlajoletSedgewick}) is the solution of
$\mathsf Y=ze^{\mathsf Y}$ taking values in $[0,1]$ for $0\le z\le e^{-1}$.
Moreover, $z\mathsf U'(z)=\mathsf Y(z)$, so differentiation of $\mathsf Y=ze^{\mathsf Y}$, followed
by integration from zero, gives $\mathsf U=\mathsf Y-\mathsf Y^2/2$. Consequently,
\[
 c(t)=\frac{\mathsf U(te^{-t})}{t}
     =\frac{\mathsf Y(te^{-t})-\mathsf Y(te^{-t})^2/2}{t}.
\]
For $0<t\le1$, the principal solution is $\mathsf Y(te^{-t})=t$,
which proves the first identity. For $t>1$, this solution lies
strictly between zero and one, giving $tc(t)<1/2$.
Finally, $\mathsf Y(z)=ze^{\mathsf Y(z)}\le ez$, so
$c(t)\le \mathsf Y(te^{-t})/t\le e^{1-t}$ for $t\ge1$.
\end{proof}

\begin{proof}[Proof of Theorem~\ref{thm:power}]
Fix $S>0$, and put $\alpha=1/40$, $a=d^{-\alpha}$,
$L=\log d$, and $\vartheta_d=A(K,S+3)L^2$, with $A$ from
Lemma~\ref{lem:connectivity}. For large $d$, $1<\vartheta_d\le d^\alpha$.
Lemma~\ref{lem:exact-subcritical} gives $0\le c(t)\le1$, $|c(t)-1|\le t$ for $0<t<1$, and $c(t)\le e^{1-t}$ for $t\ge1$.

We first bound the small-$t$ integral uniformly as $s\downarrow0$.
Write $p_d(t)=\Prob(Q_{12}\le t/d)$. Each edge reduces the number of components by at most one, so
\[
 0\le1-1/d-f_d(t)\le\frac{d-1}{2}p_d(t).
\]
The conditional row law implies
\[
 p_d(t)\le\E(t/d)^{\beta_R},\qquad
 \beta_R=(d-2)\sigma_d^2/R^2.
\]
Under \eqref{eq:bounded}, $\beta_R\ge b_K>0$ deterministically
for all sufficiently large $d$. Put $\delta=d^{-1/8}$.
Except on an event of probability
$\rho_d\le C_K e^{-c_Kd^{3/4}}$, radius concentration gives
$\beta_R\ge1-C_K\delta$. Therefore, uniformly for $0<s\le S$,
\[
 \begin{split}
 s\int_0^a(1-1/d-f_d(t))t^{s-1}\,dt
 &\le\frac{ds}{2}\E\frac{d^{-\beta_R}a^{s+\beta_R}}
                                      {s+\beta_R}\\
 &\le C_K\frac{s}{s+1}a^{s+1}+C_Kd\rho_d.
 \end{split}
\]
Here $d^{1-\beta_R}a^{\beta_R}\le C_Ka$ on the good event,
because $\delta L=o(1)$. Consequently,
\[
 s\int_0^a|f_d(t)-c(t)|t^{s-1}\,dt
 \le\frac{a^s}{d}+C_K\frac{s}{s+1}a^{s+1}
                          +C_Kd\rho_d.
\]
The term $a^s/d$ accounts for the different limits
$f_d(0+)=1-1/d$ and $c(0+)=1$.

On $[a,\vartheta_d]$, Theorem~\ref{thm:components} gives
\[
 s\int_a^{\vartheta_d}|f_d(t)-c(t)|t^{s-1}\,dt
 \le C_Kd^{-\alpha}(\vartheta_d^s-a^s).
\]
For $t\ge\vartheta_d$, monotonicity of the threshold graph and
Lemma~\ref{lem:connectivity} imply $f_d(t)\le d^{-(S+3)}$.
Hence, uniformly for $0<s\le S$,
\[
 s\int_{\vartheta_d}^d f_d(t)t^{s-1}\,dt\le d^{-3},\qquad
 s\int_{\vartheta_d}^\infty c(t)t^{s-1}\,dt
          \le C_S e^{-\vartheta_d/2}.
\]
Combining these estimates with \eqref{eq:power-kruskal} yields
\[
 |\mathcal C_d(s)-\mathcal C(s)|\le C_{K,S}d^{-\alpha}\vartheta_d^s
 \le C_{K,S}d^{-\alpha}L^{2s},
\]
as asserted. For iid uniform edge weights, $p_d(t)=t/d$,
and the same component and connectivity estimates apply, so the
proof is unchanged.
\end{proof}

\begin{proposition}[Laplace identity and a tail bound]\label{prop:laplace}
For every $d$,
\begin{equation}\label{eq:laplace-exact}
 F_d=\sqrt{\frac L\pi}\int_0^\infty
       s^{-1/2}e^{-Ls}\mathcal C_d(s)\,ds.
\end{equation}
Under \eqref{eq:bounded},
\begin{equation}\label{eq:laplace-large-s}
 \sqrt{\frac L\pi}\int_{1/2}^\infty
       s^{-1/2}e^{-Ls}\mathcal C_d(s)\,ds
 =O_K(d^{-1/2}\sqrt L).
\end{equation}
\end{proposition}

\begin{proof}
For $x>0$,
\[
 x^{-1/2}=\pi^{-1/2}\int_0^\infty s^{-1/2}e^{-sx}\,ds.
\]
Tonelli's theorem and the identity
\[
 d^{-s}\mathcal C_d(s)=d^{-1}\E\sum_{e\in\mathcal T_d}e^{-sX_e}
\]
prove \eqref{eq:laplace-exact}.

For the large-$s$ part, put $\vartheta_d=A(K,6)L^2$ and
$\mathcal E_d=\{\mathcal G_d(\vartheta_d)\text{ is connected}\}$.
Lemma~\ref{lem:connectivity} gives
$\Prob(\mathcal E_d^c)\le d^{-6}$.
On $\mathcal E_d$, Kruskal's algorithm selects every MST edge
by threshold $\vartheta_d$, so
\[
 X_e\ge \ell_d:=\log(d/\vartheta_d),\qquad e\in\mathcal T_d.
\]
For $x>0$,
\[
 I(x):=\int_{1/2}^\infty s^{-1/2}e^{-sx}\,ds
 \le\frac{\sqrt2}{x}e^{-x/2}.
\]
Thus the contribution on $\mathcal E_d$ to the integral in
\eqref{eq:laplace-large-s} is at most
\[
 C\sqrt L\,\ell_d^{-1}e^{-\ell_d/2}
 \le C_Kd^{-1/2}\sqrt L,
\]
since $\ell_d\sim L$ and $e^{-\ell_d/2}=\sqrt{\vartheta_d/d}$.

To control the exceptional event, let
$Y_d=\sum_{j=2}^d X_{1j}^{-1/2}$ be the normalized reciprocal weight of a fixed
star. The conditional beta distribution of a single edge and
\eqref{eq:inverse-radius-general} give, with $n=d-1$,
\[
 \E X_{12}^{-3/4}
 =\sigma_d^{3/2}\E R^{-3/2}
       \frac{\Gamma(1/4)\Gamma(n)}{\Gamma(n-3/4)}\le C_K.
\]
In the iid model, the same bound follows from $\E X_{12}^{-3/4}=\Gamma(1/4)$. Minkowski's inequality therefore yields $\|Y_d\|_{3/2}\le C_Kd$;
no independence between the star edges is needed.
The function $Q\mapsto I(-\log Q)$ is increasing, so its MST
edge set is again $\mathcal T_d$. Since
$I(x)\le\sqrt\pi\,x^{-1/2}$, comparison with the fixed star gives
\[
 \sum_{e\in\mathcal T_d}I(X_e)\le\sqrt\pi\,Y_d.
\]
H\"older's inequality bounds the exceptional contribution by
\[
 \frac{\sqrt L}{d}\E[Y_d\ind_{\mathcal E_d^c}]
 \le C_K\sqrt L\,\Prob(\mathcal E_d^c)^{1/3}
 \le C_K\sqrt L\,d^{-2}.
\]
This proves \eqref{eq:laplace-large-s}.
\end{proof}

\begin{proof}[Proof of Theorem~\ref{thm:universality}]
Use Proposition~\ref{prop:laplace} and Theorem~\ref{thm:power}.
For $0<s\le1/2$, put $\alpha=1/40$. The uniform bound
\eqref{eq:power-main} gives
\[
 \begin{split}
 &\sqrt L\int_0^{1/2}s^{-1/2}e^{-Ls}
                    |\mathcal C_d(s)-\mathcal C(s)|\,ds\\
 &\qquad\le C_Kd^{-\alpha}\sqrt L
       \int_0^{1/2}s^{-1/2}e^{-(L-2\log L)s}\,ds
 \le C_Kd^{-\alpha}.
 \end{split}
\]
To identify the reference profile, define for all $s>0$
\[
 C_{L}^{\rm cut}(s)=s\int_0^{e^{L/2}}c(t)t^{s-1}\,dt.
\]
Tonelli and $\Gamma(3/2)=\sqrt\pi/2$ give the exact identity
\[
 \sqrt{\frac L\pi}\int_0^\infty
 s^{-1/2}e^{-Ls}C_L^{\rm cut}(s)\,ds
 =\frac{\sqrt L}{2}\int_0^{e^{L/2}}
      \frac{c(t)}{(L-\log t)^{3/2}}\frac{dt}{t}
 =\Phi(L).
\]
For $0<s\le1/2$, replacing $C_L^{\rm cut}(s)$ by $\mathcal C(s)$ has
exponentially small error, since their difference is at most
$Cs e^{-e^{L/2}}$. The part with $s\ge1/2$ is
$O(d^{-1/2}/\sqrt L)$: for $y\ge1$,
\[
 \int_{1/2}^\infty s^{1/2}e^{-sy}\,ds
 \le C e^{-y/2}/y,
\]
and here $y=L-\log t\ge L/2$. Thus that part is bounded by
\[
 \frac{C}{\sqrt L}d^{-1/2}
       \int_0^{\sqrt d}c(t)t^{-1/2}\,dt
 \le \frac{C}{\sqrt L}d^{-1/2}.
\]
Combining these estimates proves \eqref{eq:integral-main}. Applying the same argument to the iid model gives \eqref{eq:universality}.
\end{proof}

\paragraph*{The range of the limiting transform}
The untruncated expression
$\int_0^\infty s^{-1/2}e^{-Ls}\mathcal C(s)\,ds$ diverges for every fixed
$L$: already the $k=1$ term gives $\mathcal C(s)\ge\Gamma(s+1)$. Accordingly, $\mathcal C(s)$ is integrated on a bounded $s$ interval,
or its component integral is truncated in $t$ as above. The exact
identity at finite matrix size integrates $\mathcal C_d(s)$ over all $s>0$.

\paragraph*{The finite-size centering term}

With $\kappa_d(t)$ as in Section~\ref{sec:component-universality}, define
\[
 h_d(t)=\left(\frac{L}{L-\log t}\right)^{3/2},
\]
and the centered component density
\begin{equation}\label{eq:centered-density}
 g_d(t)=\frac{\E\kappa_d(t)-1}{d}
       -\frac{d-1}{d}\ind_{\{t<1\}}.
\end{equation}
This subtraction makes $g_d(t)\to0$ as $t\downarrow0$.

\begin{lemma}[Centered component-count identity]\label{lem:centered-kruskal}
One has
\begin{equation}\label{eq:centered-kruskal}
 F_d=\frac{d-1}{d}
       +\frac1{2L}\int_0^d g_d(t)h_d(t)\,\frac{dt}{t}.
\end{equation}
\end{lemma}
\begin{proof}
For the graph retaining weights at most $w$, let $\widehat\kappa_d(w)$ be
its component count. The component-count identity
$M_d=\int_0^\infty(\widehat\kappa_d(w)-1)\,dw$, used by
\citet{AvramBertsimas} and \citet{Janson}, and the change of variables
\[
 w=w_d(t)=\frac1{\sigma_d\sqrt{L-\log t}},\qquad
 dw=\frac{dt}{2\sigma_dt(L-\log t)^{3/2}}
\]
give $\widehat\kappa_d(w_d(t))=\kappa_d(t)$ away from the finitely
many edge thresholds and yield the formula before centering.
The subtracted term integrates
exactly, since
\[
 \int_0^1 h_d(t)\,\frac{dt}{t}=2L.
\]
All integrands before centering are nonnegative, so the change of
variables is justified by Tonelli's theorem. The inverse-edge moment
bound ensures $\E M_d<\infty$.
\end{proof}

The centering term $(d-1)/d$ explains one contribution of order $1/d$.
Other corrections at that order are not determined here.

\section{The reciprocal weights and their asymptotic expansion}
\label{sec:coefficients}

Write $h(t)=c(t)-\ind_{\{t<1\}}$. We evaluate the profile
\eqref{eq:Phi} and identify the Mellin singularity responsible for
the divergence of its asymptotic expansion.

Lemma~\ref{lem:exact-subcritical} gives a stable alternative formula for
\eqref{eq:Phi}:
\begin{equation}\label{eq:phi-stable}
 \Phi(L)=\frac12+\frac{\sqrt{\pi L}}2
             \operatorname{erfcx}(\sqrt L)
 +\frac1{2L}\int_1^{e^{L/2}}
 c(t)\left(1-\frac{\log t}{L}\right)^{-3/2}\frac{dt}{t},
\end{equation}
where
$\operatorname{erfc}(x)=\frac2{\sqrt\pi}\int_x^\infty e^{-u^2}\,du$
is the complementary error function and
$\operatorname{erfcx}(x)=e^{x^2}\operatorname{erfc}(x)$.
Indeed, the part below one is
\[
 -\frac1{4L}\int_0^\infty
    e^{-x}(1+x/L)^{-3/2}\,dx
 =-\frac12+\frac{\sqrt{\pi L}}2
                    \operatorname{erfcx}(\sqrt L).
\]
Thus this part can be evaluated without numerical integration near
$t=0$ or evaluation of the product $e^L\operatorname{erfc}(\sqrt L)$.

\begin{proposition}[Mellin transform and coefficient evaluation]
\label{prop:coefficient-evaluation}
For $\Re s>-1$, define
\begin{equation}\label{eq:mellin-component}
 F(s)=\int_0^\infty h(t)t^{s-1}\,dt.
\end{equation}
Then
\begin{align}
 F(s)&=\Gamma(s)-\frac1s+
       \sum_{k\ge2}\frac{\Gamma(k+s-1)}{k!\,k^{s+1}},
       \label{eq:mellin-series}\\
 F(s)&=-\frac1{2(s+1)}+G(s),\qquad
 G(s)=\int_1^\infty c(t)t^{s-1}\,dt.
       \label{eq:mellin-endpoint}
\end{align}
The singularity at zero in \eqref{eq:mellin-series} is removable,
and $G$ is entire. Thus $F$ extends meromorphically to the whole
plane with exactly one pole, at $s=-1$, of residue $-1/2$.
The logarithmic moments are $\mu_r=F^{(r)}(0)$ and
\begin{equation}\label{eq:universal-coefficients}
 a_{r+1}=\frac{(3/2)_r}{2r!}\mu_r,
 \qquad
 a_1=1-\frac\gamma2-\frac{\pi^2}{12}
     \approx-0.1111.
\end{equation}
\end{proposition}
\begin{proof}
The $k=1$ summand of $c$ is $e^{-t}$, whose regularized Mellin
integral is $\Gamma(s)-1/s$. For $k\ge2$, termwise integration gives
\[
 \frac{k^{k-2}}{k!}\int_0^\infty
 t^{k+s-2}e^{-kt}\,dt
 =\frac{\Gamma(k+s-1)}{k!\,k^{s+1}}.
\]
The resulting series is locally uniformly convergent in
$\Re s>-1$, since its summands are $O(k^{-3})$ uniformly on
compact sets. This proves \eqref{eq:mellin-series}.
Lemma \ref{lem:exact-subcritical} gives
\eqref{eq:mellin-endpoint} directly. Its exponential upper bound
makes $G$ entire, by differentiation under the integral on every
compact subset of the complex plane.

At zero, the regularized gamma term is $-\gamma$, and
\[
 \sum_{k\ge2}\frac1{k^2(k-1)}
 =\sum_{k\ge2}\left(\frac1{k-1}-\frac1k-\frac1{k^2}\right)
 =2-\zeta(2).
\]
Hence $F(0)=2-\gamma-\zeta(2)$, proving the value of $a_1$.
\end{proof}

\begin{proposition}[Asymptotic but divergent series]
\label{prop:coefficient-divergence}
For each fixed integer $N\ge1$,
\begin{equation}\label{eq:phi-expansion}
 \Phi(L)=1+\sum_{r=0}^{N-1}\frac{a_{r+1}}{L^{r+1}}
                 +O_N(L^{-N-1}).
\end{equation}
Every coefficient has the alternating sign
$\operatorname{sgn}(a_{r+1})=(-1)^{r+1}$. More precisely, for
every fixed $\rho>0$,
\begin{equation}\label{eq:large-order-coefficients}
 a_{r+1}=\frac{(-1)^{r+1}}4(3/2)_r
          \left[1+O_\rho(\rho^{-r})\right],\qquad r\to\infty.
\end{equation}
In particular the formal inverse-logarithmic series has zero radius
of convergence; it is a Poincar\'e asymptotic expansion.
\end{proposition}
\begin{proof}
On the part $t<1$, write $x=-\log t$. Repeated integration by
parts gives the exact remainder identity
\begin{align*}
 \int_0^\infty e^{-x}(1+x/L)^{-3/2}\,dx
 &=\sum_{r=0}^{N-1}\frac{(-1)^r(3/2)_r}{L^r}\\
 &\quad+\frac{(-1)^N(3/2)_N}{L^N}
       \int_0^\infty e^{-x}(1+x/L)^{-N-3/2}\,dx.
\end{align*}
The last integral lies between zero and one. On $1\le t\le e^{L/2}$,
the quantity $x=(\log t)/L$ lies in $[0,1/2]$, so Taylor's theorem
gives
\[
 (1-x)^{-3/2}
 =\sum_{r=0}^{N-1}\frac{(3/2)_r}{r!}x^r+O_N(x^N).
\]
All positive logarithmic moments of $c(t)\,dt/t$ are finite.
The integrated Taylor remainder is therefore $O_N(L^{-N-1})$;
extending each polynomial moment from $e^{L/2}$ to infinity introduces an error
less than every inverse power of $e^L$. This proves
\eqref{eq:phi-expansion}.

For the coefficient signs and large-order behavior,
\eqref{eq:mellin-endpoint} gives exactly
\begin{equation}\label{eq:moment-splitting}
 \mu_r=\frac{(-1)^{r+1}r!}{2}
       +\int_1^\infty c(t)(\log t)^r\,\frac{dt}{t}.
\end{equation}
The second term is strictly positive and, by
$c(t)<1/(2t)$ for $t>1$, strictly smaller than
\[
 \frac12\int_1^\infty(\log t)^r\frac{dt}{t^2}
 =\frac{r!}{2}.
\]
This proves the alternating signs for every $r$. Since $G$ is
entire, Cauchy's estimate on $|s|=\rho$ gives
$G^{(r)}(0)=O_\rho(r!\rho^{-r})$. Substitution into
\eqref{eq:moment-splitting} proves
\eqref{eq:large-order-coefficients}. In particular,
$a_{r+2}/a_{r+1}\sim-(r+3/2)$, so the series terms fail to tend
to zero for any fixed finite $L$.
\end{proof}

\subsection{The Avram--Bertsimas constant and Euclidean dimension}\label{sec:geometric}
The function $\mathcal C$ also describes MST weights for other
independent edge distributions. \citet[Theorem~8 and the following remark]{AvramBertsimas} consider bounded independent edge lengths with
distribution function $G(x)\sim b_Dx^D$ at the origin. Let
$L_n^{\rm ind}$ denote the MST weight on $n$ vertices for this
independent model. Their formula is
\[
 \lim_{n\to\infty}\frac{\E L_n^{\rm ind}}{n^{1-1/D}}
 =b_D^{-1/D}\,\mathcal C(1/D),\qquad
 \mathcal C(s)=s\sum_{k\ge1}
       \frac{\Gamma(k+s-1)}{k!\,k^{s+1}}
       =s\int_0^\infty c(t)t^{s-1}\,dt.
\]
Their power-weighted version gives the same function at $s=a/D$,
$0<a\le D$, with prefactor $b_D^{-a/D}$.
Avram and Bertsimas also attribute the independent-model series to the
earlier analysis of Prim's algorithm by \citet{Timofeev}.
For the regularized Mellin transform used here,
\[
 \mathcal C(s)=1+sF(s),\qquad
 \mathcal C(1)=\zeta(3),\qquad
 \mathcal C(0)=1,\qquad
 \mathcal C'(0)=F(0)=2-\gamma-\zeta(2).
\]
The reciprocal-weight expansion uses the behavior of this classical
moment function at $s=0$. Theorem~\ref{thm:power} supplies the uniform
comparison needed to transfer it to Haar-conjugated matrices.

Since $s=1/D$, letting $s\downarrow0$ gives the behavior of
the independent-model constant as $D\to\infty$. In particular,
\[
 \mathcal C(1/D)=1+\frac{2-\gamma-\zeta(2)}D+O(D^{-2}).
\]
For comparison with the geometric constant in
Proposition~\ref{cor:geometric-endpoint}, let $\beta_D^{\rm ind}$ be the
independent-model constant with lower-tail coefficient $\omega_D$.
The normalization is
$\omega_D^{1/D}\beta_D^{\rm ind}=\mathcal C(1/D)$.
\citet[Corollary~1 and Proposition~2]{Penrose} proves that
the geometric component density $c_D^{\rm geom}(t)$ satisfies
\[
 c_D^{\rm geom}(t)\ge c(t),\qquad
 c_D^{\rm geom}(t)\longrightarrow c(t)
 \quad(D\to\infty)
\]
at every fixed mean degree $t>0$. Here $c_D^{\rm geom}(t)$ is the
expected reciprocal size of the origin cluster in the Poisson
geometric graph of intensity $t/\omega_D$ and connection radius one,
with an extra point at the origin; an infinite cluster contributes
zero. His Theorem~2 restates the
Avram--Bertsimas component integrals and yields
$\beta_D^{\rm geom}\ge\beta_D^{\rm ind}$.
The leading relation $\beta_D^{\rm geom}/\beta_D^{\rm ind}\to1$
already follows from the bounds of \citet{BertsimasVanRyzin}.
There is a more precise edge-law interpretation in
\citet[Proposition~3 and Theorem~2, (7)--(8)]{PenroseHighDim}.
Let $\tau_D$ have the limiting empirical distribution of
$n\omega_D|e|^D$ over Euclidean MST edges, as $n\to\infty$ at fixed $D$.
Penrose proves
\[
 \Prob(\tau_D\le t)\longrightarrow 1-c(t),\qquad t>0.
\]
His Lemma~3 and equation~(26) in the proof of Lemma~4 imply
$\Prob(\tau_D>t)\le C_*e^{-c_*t}$ for $t\ge t_0$ and $D\ge D_0$,
with positive constants independent of $D$. Counting Poisson neighbours
also gives $\Prob(\tau_D\le t)\le t/2$. These two bounds justify
convergence of logarithmic moments and the Taylor expansion of
$\E\tau_D^{1/D}=\omega_D^{1/D}\beta_D^{\rm geom}$ at exponent zero.
Thus Proposition~\ref{cor:geometric-endpoint} is an explicit consequence
of Penrose's edge-law limit and tail estimates: its coefficient is
$\E\log\tau=2-\gamma-\zeta(2)$. The component-count proof below
obtains the same conclusion with a simpler integrable bound.

The cube and flat unit-torus MST constants agree for each fixed
$D\ge2$, by \citet[Theorem~1]{Jaillet}. Thus the torus
constant used by Avram--Bertsimas is the cube constant in
\citet[Theorem~2]{Penrose} and in Proposition~\ref{cor:geometric-endpoint}.

\begin{proof}[Proof of Proposition~\ref{cor:geometric-endpoint}]
Use the continuum-percolation parametrization of
\citet[Sections~1--3]{Penrose}. Let $\mathcal P_{D,t}$ be a
homogeneous Poisson process in $\mathbb R^D$ of intensity $t/\omega_D$,
with one extra point inserted at the origin. Join points whose
Euclidean distance is less than one, and let $\mathcal K_{D,t}$ be
the component containing the origin. Define
$f_k^{(D)}(t)=\Prob(|\mathcal K_{D,t}|=k)$ for $k=1,2,\ldots$.
Then the density introduced above is
\[
 c_D^{\rm geom}(t)=\E\frac1{|\mathcal K_{D,t}|}
       =\sum_{k\ge1}\frac1k f_k^{(D)}(t),
 \qquad \frac1\infty:=0.
\]
Penrose's Corollary~1, equation~(7), with his parameter $K=1$
(unrelated to the constant in \eqref{eq:bounded}), and
Proposition~2, equation~(8), give respectively
\[
 c_D^{\rm geom}(t)\longrightarrow c(t)\quad(D\to\infty),\qquad
 c_D^{\rm geom}(t)\ge c(t),\qquad t>0.
\]

There is also a dimension-independent upper bound. The number $N_t$
of Poisson points at distance less than one from the origin is
Poisson with mean $t$, since the intensity is $t/\omega_D$. Every one of
these points belongs to the origin component, so, pathwise,
$|\mathcal K_{D,t}|\ge N_t+1$. This remains valid when the component
is infinite, with the convention above. Therefore
\begin{equation}\label{eq:geometric-degree-profile}
 c(t)\le c_D^{\rm geom}(t)\le\E\frac1{N_t+1}
 =e^{-t}\sum_{m=0}^{\infty}\frac{t^m}{m!(m+1)}
 =\frac{1-e^{-t}}{t}=:r(t).
\end{equation}

The Euclidean MST representation restated by
\citet[Theorem~2, equation~(12)]{Penrose} is
\[
 \omega_D^{1/D}\beta_D^{\rm geom}
 =\frac1D\int_0^\infty c_D^{\rm geom}(t)t^{1/D-1}\,dt.
\]
Subtracting the corresponding independent-model integral yields
\begin{equation}\label{eq:geometric-profile-difference}
 D\left[\omega_D^{1/D}\beta_D^{\rm geom}-\mathcal C(1/D)\right]
 =\int_0^\infty[c_D^{\rm geom}(t)-c(t)]t^{1/D-1}\,dt.
\end{equation}
For $0<t\le1$, the identity $c(t)=1-t/2$ gives
\[
 0\le r(t)-c(t)
 =\frac{1-e^{-t}}t-1+\frac t2\le\frac{t^2}{6}.
\]
For $t\ge1$, one has $0\le r(t)-c(t)\le1/t$. Hence, for every
$D\ge2$, the nonnegative integrand in
\eqref{eq:geometric-profile-difference} is bounded by
\[
 \frac t6\,\ind_{\{0<t\le1\}}
 +t^{-3/2}\,\ind_{\{t>1\}}.
\]
This is integrable and independent of $D$, while pointwise convergence
of $c_D^{\rm geom}$ makes the integrand tend to zero. Dominated convergence
therefore proves the first equality in
\eqref{eq:penrose-degree-correction}. The second follows from
$\mathcal C(s)=1+(2-\gamma-\zeta(2))s+O(s^2)$ at zero.

\end{proof}

The order of limits matters: the number of sampled points first
tends to infinity at fixed spatial dimension, and then $D\to\infty$.
The argument gives an explicit coefficient but no rate for its remainder.

\subsection{Watson's lemma and the reciprocal-weight expansion}
For fixed $S>0$, define the truncated transform
\begin{equation}\label{eq:truncated-watson}
 \Psi_S(L)=\sqrt{\frac L\pi}\int_0^S
       s^{-1/2}e^{-Ls}\mathcal C(s)\,ds.
\end{equation}
This is an application of Watson's lemma; see
\citet[pp.~71--72]{Olver}. Here Taylor's theorem at zero and
integration of finitely many monomials give a direct proof. If $[s^n]\mathcal C(s)$ denotes
the Taylor coefficient at zero, the coefficient of $L^{-n}$ is
$\Gamma(n+1/2)[s^n]\mathcal C(s)/\sqrt\pi$.
Consequently,
\begin{equation}\label{eq:watson-coefficients}
 a_{r+1}=\frac{(3/2)_r}{2}[s^{r+1}]\mathcal C(s)
         =\frac{(3/2)_r\mu_r}{2r!},\qquad
 a_1=\frac{\mathcal C'(0)}2.
\end{equation}
For a fixed truncation after degree $N$, the Taylor remainder is
$O_{N,S}(s^{N+1})$, whose normalized Laplace integral is
$O_{N,S}(L^{-N-1})$. Extending the finitely many polynomial integrals
to infinity contributes an exponentially small error.
Thus the logarithmic coefficients are the Taylor coefficients of
the classical moment function, multiplied by the gamma factors
from integrating against the weight $s^{-1/2}e^{-Ls}$.

Equation~\eqref{eq:mellin-endpoint} yields
\[
 \mathcal C(s)=\frac12+\frac1{2(s+1)}+sG(s).
\]
Its pole at $s=-1$ has residue $+1/2$, giving the factorial growth
in \eqref{eq:large-order-coefficients} after Watson's transformation.

\begin{proof}[Proof of Corollary~\ref{thm:expansion}]
Proposition~\ref{prop:coefficient-divergence} gives the expansion
of $\Phi$ with coefficients evaluated in
Proposition~\ref{prop:coefficient-evaluation}.
Theorem~\ref{thm:universality} transfers it to both random-matrix
models because $d^{-1/40}=o((\log d)^{-m})$ for every fixed $m$.
The coefficient asymptotic in Proposition~\ref{prop:coefficient-divergence}
also proves the asserted divergence.
\end{proof}

For numerical evaluation of the second coefficient,
\begin{align}
 \mu_1&=\frac{\gamma^2}{2}+\frac{\pi^2}{12}
  +\sum_{k\ge2}\frac{\psi(k-1)-\log k}{k^2(k-1)},\label{eq:mu-one}\\
 a_2&=\frac34\mu_1\approx0.4627,\nonumber
\end{align}
where $\psi=\Gamma'/\Gamma$. Differentiating
\eqref{eq:mellin-series} gives this absolutely convergent expression.

Appendix~\ref{sec:numerics} gives a deterministic evaluation of the
reciprocal-weight profile. The recursive-tree expansion is given in
Appendix~\ref{app:recursive}.

\section{Spectral growth and vanishing MST weight}\label{sec:spectral-growth}
For fixed $a>0$, the spectrum $\lambda_{k,d}=k^a$ satisfies
\eqref{eq:bounded}, and Riemann sums give
\begin{equation}\label{eq:power-variance}
 v_d\sim\frac{a^2d^{2a}}{(2a+1)(a+1)^2}.
\end{equation}
Consequently,
\begin{equation}\label{eq:power-leading}
 \E M_d\sim\frac{(a+1)\sqrt{2a+1}}{a}
                 \frac{d^{3/2-a}}{\sqrt{\log d}}.
\end{equation}
For $a=3/2$, the Riemann-sum error in $v_d$ is relative $O(d^{-1})$,
so Theorems \ref{thm:universality} and \ref{thm:main} give
\begin{equation}\label{eq:critical-second}
 \E M_d=\frac{10}{3\sqrt L}
       \left(\Phi(L)+O(d^{-1/40})\right),
 \qquad
 \sqrt L\,M_d\longrightarrow\frac{10}{3}\quad\text{in }L^1.
\end{equation}
Under \eqref{eq:bounded}, Theorem \ref{thm:main} also gives the
criterion for the MST weight to tend to zero:
\begin{equation}\label{eq:certificate-criterion}
 M_d\longrightarrow0\text{ in probability (equivalently in }L^1\text{)}
 \quad\Longleftrightarrow\quad
 \frac{v_d\log d}{d^3}\longrightarrow\infty.
\end{equation}
For necessity, if the deterministic scale
$d^{3/2}/\sqrt{v_d\log d}$ fails to tend to zero, it has a subsequence
bounded below; \eqref{eq:main} then rules out convergence of $M_d$ to
zero in probability. Sufficiency follows directly from \eqref{eq:main}.

In particular, for fixed $b\in\mathbb R$ and
$\lambda_k=k^{3/2}/(\log k)^b$, $k\ge2$ (with any fixed $\lambda_1$),
regular-variation Riemann sums give
\[
 v_d\sim\frac9{100}\frac{d^3}{(\log d)^{2b}},\qquad
 \E M_d\sim\frac{10}{3}(\log d)^{b-1/2}.
\]
Thus $M_d\to0$ in $L^1$ for $b<1/2$. At $b=1/2$,
$M_d\to10/3$ in $L^1$; for $b>1/2$ its deterministic first-order
scale diverges.

\section{Further questions}\label{sec:open}
A natural next question is whether the Haar threshold graph has the
same connectivity scale as the independent graph. Under
\eqref{eq:bounded}, does every fixed $\varepsilon>0$ satisfy
\[
 \Prob\bigl(\mathcal G_d((1+\varepsilon)\log d)
                 \text{ is connected}\bigr)\longrightarrow1?
\]
Lemma~\ref{lem:connectivity} gives connectivity at threshold
parameter $A_K(\log d)^2$. The additional logarithm in that
proof comes from the contribution of a single exceedance to the
squared function $F_S$ in the proof of Lemma~\ref{lem:connectivity}.

Theorem~\ref{thm:power} concerns expectations. A corresponding law
of large numbers for the Haar MST with weights $Q_e$ would strengthen
the result at $s=1$. Extending the comparison to real orthogonal
conjugation is another natural question; the one-row distribution
then has Dirichlet parameters $1/2$ instead of $1$.

Sharp finite-size corrections are known in the iid uniform model:
\citet{CooperEtAl} obtain an expansion through order $d^{-4/3}$.
For Haar dependence and for reciprocal weights, the comparisons
in Appendix~\ref{sec:numerics} suggest corrections much smaller than
our bound, but proving them requires more precise component estimates.
\citet[Table~3]{CM} report different fitted variance slopes for the
two models at the simulated sizes. Whether these reflect distinct
asymptotic fluctuations remains open. A variance comparison and a
central limit theorem require control of joint component counts
across thresholds, beyond the expectation estimates proved here.

\appendix

\section{The recursive-tree expansion}\label{app:recursive}
For the recursive tree, the limiting profile is the function
$r(t)=(1-e^{-t})/t$ already used in \eqref{eq:geometric-degree-profile}.
Its regularized Mellin transform is
\begin{equation}\label{eq:mellin-recursive}
 \mathcal R(s):=\int_0^\infty
 [r(t)-\mathbf1_{\{t<1\}}]t^{s-1}\,dt
 =\frac{\Gamma(s)}{1-s}-\frac1s,
 \qquad -1<\Re s<1.
\end{equation}
Indeed, integration by parts gives
$\int_0^\infty(1-e^{-t})t^{s-2}\,dt=\Gamma(s)/(1-s)$ for
$0<\Re s<1$, and subtraction of $1/s$ extends the formula to the
displayed strip. The same expansion of $(1-\log t/L)^{-3/2}$ therefore gives
\begin{equation}\label{eq:recursive-expansion}
 \frac{\sigma_d\sqrt{\log d}}d\E M_d^{\rm rec}
 =1+\frac{a_1^{\rm rec}}{\log d}+\frac{a_2^{\rm rec}}{(\log d)^2}
   +O_K((\log d)^{-3}),
\end{equation}
where
\begin{align*}
 a_1^{\rm rec}&=\frac{1-\gamma}{2}\approx0.2114,\\
 a_2^{\rm rec}&=\frac34\left(1-\gamma+\frac{\gamma^2}{2}
                              +\frac{\pi^2}{12}\right)
      \approx1.0589.
\end{align*}
In particular,
\begin{equation}\label{eq:tree-gap}
 a_1^{\rm rec}-a_1=\frac{\pi^2}{12}-\frac12
 \approx0.3225>0.
\end{equation}

Here is a direct justification of the recursive expansion, independent
of component-count comparison for the MST. For the Gaussian recursive
tree weight $M_d^{\rm rec,iid}$, let $r_d(t)$ be $1/d$ times the expected
number of its edges with transformed weight $Q_e>t/d$. Its exact
survival profile is
\[
 r_d(t)=\frac1d\sum_{m=1}^{d-1}(1-t/d)^m
       =\frac{(1-t/d)-(1-t/d)^d}{t},\qquad 0<t<d.
\]
The tail integral for $M_d^{\rm rec,iid}$ uses the same weight $h_d(t)$ as
\eqref{eq:centered-kruskal}, with $r_d$ in place of
$f_d(t)=(\E\kappa_d(t)-1)/d$. Subtract
$(1-1/d)\mathbf1_{\{t<1\}}$ before expanding.
On $0<t<1$, the difference of this centered profile from
$r(t)-1$ is $O(t/d)$; on $1\le t\le d/2$, $r_d(t)-r(t)=O(1/d)$.
These estimates follow from
\[
 r_d(t)-r(t)=-\frac1d+
             \frac{e^{-t}-(1-t/d)^d}{t}.
\]
Also $r_d(t)\le1/t$, and near $t=d$ it satisfies
$r_d(t)\le(1-t/d)/t$. Splitting at $d^{-\alpha}$, $d^\alpha$,
and $d/2$, with any fixed $0<\alpha<1$, makes the two extreme tails
polynomially small in $d$, even after multiplication by any fixed
power of $\log d$. On the central interval the Taylor remainder is
controlled by the finite logarithmic moments of
$r(t)-\mathbf1_{\{t<1\}}$. This proves an expansion to every fixed
order, with coefficients obtained by replacing $F$ by $\mathcal R$ in
\eqref{eq:universal-coefficients}.

For the Haar recursive tree, the exact one-row identity gives
\[
 \frac{\E M_d^{\rm rec}}{\E M_d^{\rm rec,iid}}
 =\sigma_d\E R_1^{-1}
   \frac{\Gamma(d-1)}{\Gamma(d-3/2)}
 =1+O_K(d^{-1/4}).
\]
For the last estimate, use radius concentration with relative tolerance
$d^{-1/4}$, and control the exceptional contribution by H\"older's
inequality and the bounded $3/2$ inverse moment from
\eqref{eq:inverse-radius-general}. The gamma ratio and the exact
formula $\sigma_d^2=dv_d/(d^2-1)$ contribute only $O(d^{-1})$.
The polynomial error is smaller than every fixed inverse power of
$\log d$, establishing \eqref{eq:recursive-expansion} for the Haar model.

Combining the two expansions also gives
\[
 \frac{\E M_d^{\rm rec}}{\E M_d}-1
 =\frac{a_1^{\rm rec}-a_1}{\log d}
  +\frac{a_2^{\rm rec}-a_2+a_1^2-a_1a_1^{\rm rec}}{(\log d)^2}
  +O_K((\log d)^{-3}).
\]
The two coefficients are approximately $0.3225$ and $0.6320$.
For evaluation of the iid recursive tree at finite matrix size,
an exact alternative is
\[
 \frac{\sigma_d\sqrt L}{d}\E M_d^{\rm rec,iid}
 =\frac{\sqrt L}{2d}\int_0^\infty x^{-3/2}
    \sum_{m=1}^{d-1}(1-e^{-x})^m\,dx.
\]
This follows from the distribution of the exponential maximum $Z_m$.
The ratio expansion above describes the limiting coefficients; the
exact integral is available when values at finite matrix size are needed.

\section{Numerical evaluation and published means}\label{sec:numerics}
For the spectrum $\lambda_k=k^{3/2}$, \citet[Tables~1--2]{CM}
report Haar and iid means at sixteen matrix sizes. Let
\[
 P_d=\frac{d}{\sigma_d\sqrt{\log d}}\Phi(\log d),\qquad
 \delta_d^{\rm model}=d\left(\frac{\widehat M_d^{\rm model}}{P_d}-1\right),
\]
where $\widehat M_d^{\rm model}$ is the reported sample mean.
The exact finite-$d$ spectral variance is used in $\sigma_d$.
Table~\ref{tab:published-means} reproduces their means and compares them with
the deterministic profile. All displayed means are rounded to five
decimal places; the residuals use unrounded profile values. Parentheses
give estimated Monte Carlo standard errors, computed as
$d\,s_d/(P_d\sqrt{n_d})$ from the reported sample standard deviation
$s_d$ and number $n_d$ of simulation runs for each model and matrix size.

\begin{table}[tb]
\caption{Published means and comparison with the deterministic profile. Parentheses give estimated Monte Carlo standard errors.}
\label{tab:published-means}
\centering
\begin{tabular}{rrrrrr}
\toprule
$d$ & iid mean & Haar mean & $P_d$ & $\delta_d^{\rm iid}$ & $\delta_d^{\rm Haar}$\\
\midrule
100 & 1.52019 & 1.52018 & 1.53394 & $-0.90\;(0.03)$ & $-0.90\;(0.03)$\\
142 & 1.47117 & 1.47099 & 1.48036 & $-0.88\;(0.03)$ & $-0.90\;(0.03)$\\
200 & 1.42600 & 1.42629 & 1.43289 & $-0.96\;(0.04)$ & $-0.92\;(0.03)$\\
283 & 1.38448 & 1.38472 & 1.38900 & $-0.92\;(0.04)$ & $-0.87\;(0.04)$\\
400 & 1.34559 & 1.34589 & 1.34894 & $-0.99\;(0.05)$ & $-0.90\;(0.04)$\\
566 & 1.30997 & 1.31002 & 1.31198 & $-0.87\;(0.05)$ & $-0.84\;(0.05)$\\
800 & 1.27646 & 1.27646 & 1.27796 & $-0.94\;(0.06)$ & $-0.94\;(0.06)$\\
1131 & 1.24535 & 1.24533 & 1.24641 & $-0.96\;(0.07)$ & $-0.98\;(0.07)$\\
1600 & 1.21635 & 1.21629 & 1.21701 & $-0.87\;(0.08)$ & $-0.95\;(0.07)$\\
2263 & 1.18915 & 1.18918 & 1.18962 & $-0.89\;(0.09)$ & $-0.84\;(0.08)$\\
3200 & 1.16366 & 1.16367 & 1.16401 & $-0.96\;(0.10)$ & $-0.93\;(0.09)$\\
4525 & 1.13972 & 1.13978 & 1.13999 & $-1.06\;(0.11)$ & $-0.83\;(0.11)$\\
6400 & 1.11724 & 1.11721 & 1.11738 & $-0.83\;(0.13)$ & $-1.00\;(0.13)$\\
9051 & 1.09595 & 1.09594 & 1.09609 & $-1.13\;(0.15)$ & $-1.21\;(0.29)$\\
12800 & 1.07588 & 1.07592 & 1.07597 & $-1.04\;(0.17)$ & $-0.56\;(0.34)$\\
18101 & 1.05686 & 1.05685 & 1.05693 & $-1.11\;(0.20)$ & $-1.29\;(0.64)$\\
\bottomrule
\end{tabular}
\end{table}

The inverse-variance averages across these matrix sizes are
$-0.923\pm0.014$ for iid weights and $-0.901\pm0.013$ for Haar
weights, where the uncertainties are one estimated standard error
assuming independent simulation batches. These are descriptive averages,
not estimates of a proved limiting coefficient. The Haar sample counts
at the three largest sizes are $488$, $500$ and $200$; this explains
the larger standard errors there. The data are consistent with a
correction of order $1/d$ over the displayed range, much of which is
accounted for by the centering term in Lemma~\ref{lem:centered-kruskal}.
They do not establish that asymptotic rate, and the bound in
Theorem~\ref{thm:universality} does not determine these smaller corrections.

The effective first correction in the profile is
$A_{\mathrm{eff}}(d)=(\log d)[\Phi(\log d)-1]$.
It is approximately $-0.0565$ at $d=1131$ and $-0.0703$ at $d=18101$,
two matrix sizes used by \citet{CM}. Thus the profile helps explain
why their fitted correction near $-0.07$ differs from the limiting
coefficient $a_1\approx-0.1111$.
For example, $P_{1131}\approx1.24641$, compared with their reported
Haar mean $1.24533$.

\subsection{Quadrature check}
The truncated Laplace transform provides an independent way to
evaluate the reciprocal-weight profile. Using $S=2$ in
\eqref{eq:truncated-watson}, deterministic quadrature gives the values in
Table~\ref{tab:quadrature}.

\begin{table}[b]
\caption{Comparison of the profile and its truncated Laplace representation.}
\label{tab:quadrature}
\centering
\begin{tabular}{rrrr}
\toprule
$d$ & $\Phi(\log d)$ & $\Psi_2(\log d)$ & $|\Phi(\log d)-\Psi_2(\log d)|$\\
\midrule
1131 & 0.991969 & 0.991968 & $3.1\times10^{-7}$\\
18101 & 0.992829 & 0.992829 & $1.0\times10^{-9}$\\
\bottomrule
\end{tabular}
\end{table}

Function values are rounded to six decimal places; the differences
were computed before rounding. For stable quadrature, put $s=u^2$ to obtain
\[
 \Psi_2(L)=2\sqrt{L/\pi}\int_0^{\sqrt2}
        e^{-Lu^2}\mathcal C(u^2)\,du.
\]
The cutoff in $\Psi_2$ is necessary because the full Laplace integral of $\mathcal C$ diverges.

\section{Two consequences of the component estimates}\label{app:component-consequences}
The following observations use Theorem~\ref{thm:components} and its
proof, but are not needed for the main results. In this appendix,
$\alpha=1/40$ and $L=\log d$.

\subsection{Marginal calibration}
If the edge threshold is instead chosen to make its exact Haar
marginal probability equal to $t/d$, Theorem~\ref{thm:components} remains valid.
The one-row beta tail and radius concentration, with error
$\delta=d^{-1/8}$, show uniformly on this interval that its
normalized squared threshold differs from $\log(d/t)$ by
$O_K(\delta L)$. The discarded radius-concentration probability
is negligible relative to $p\ge d^{-1-\alpha}$. Inserting this
threshold shift into the internal and external sandwiches above
adds relative errors $O_K(k\delta L)$ and $O_K(kt\delta L)$,
respectively; both are smaller than $d^{-1/16}L$ on the stated
range. The same counting argument therefore proves
\eqref{eq:component-rate} under exact marginal calibration.

\subsection{Isolated vertices}
For fixed $t>0$, let $I_d=I_d(t)$ be the number of isolated vertices
of $\mathcal G_d(t)$. The same conditional absence estimate with
$k=1,2$ gives
\[
 \frac{\E I_d}{d}\to e^{-t},\qquad
 \frac{\E[I_d(I_d-1)]}{d(d-1)}\to e^{-2t},\qquad
 I_d/d\to e^{-t}\quad\text{in }L^2.
\]
For the second relation, apply the $k=2$ estimate to the two external
stars; requiring their shared edge to be absent changes the limiting
probability by $o(1)$.

\paragraph*{Funding.}
This work was partially supported by the Czech Science Foundation
under grant No.~24-11146S.

\section*{Statement on AI use}
OpenAI’s ChatGPT Astra 6 and Anthropic’s Claude Fable 5.2 were used during the development and revision of this work. They helped develop and verify proofs, particularly the conditional representations of Haar matrices, estimates for component events, the connectivity argument, integral representations, and asymptotic calculations. They also helped revise the exposition, write code for numerical checks, search the literature, check reference details, and edit the LaTeX source. The author retains responsibility for the mathematical claims, numerical results, source attribution, and final text.

\end{document}